%% file: JOTA_Ver2.tex
\documentclass[11pt]{article}
\usepackage[margin=.975in]{geometry}
\usepackage{times} 
\usepackage{float}
\usepackage{amsmath,amsthm,amssymb,MnSymbol}
\usepackage{graphicx}
\usepackage{hyperref} 
\usepackage[dvipsnames]{xcolor}
\hypersetup{colorlinks=true, linkcolor=Maroon, filecolor=blue, urlcolor=Brown, citecolor=Maroon} 
\usepackage[noblocks]{authblk} 
\theoremstyle{theorem}
\newtheorem{theorem}{Theorem}[section]
\newtheorem{corollary}{Corollary}[section]
\newtheorem{lemma}{Lemma}[section]

\theoremstyle{definition}
\newtheorem{definition}{Definition}[section]

\theoremstyle{remark}
\newtheorem{remark}{Remark}[section]
\newtheorem{example}{Example}[section]

\usepackage[linesnumbered,ruled]{algorithm2e}
\usepackage[export]{adjustbox}
\usepackage{subfig}
\usepackage{setspace}
\usepackage{caption} 
\usepackage{tikz}
\usetikzlibrary{arrows}
\usetikzlibrary{arrows.meta} 
\usetikzlibrary{calc,intersections,through,backgrounds}
\usetikzlibrary{positioning}
\usepackage{pgfplots}
\pgfplotsset{compat=1.16}
\usepackage{pgfplotstable}
\usepackage{filecontents}
\usepackage{lineno}

\newcommand{\N}{\mathsf N}

\newcommand{\col}[1]{\textrm{col}\{#1\}}

\begin{document}

 
%
 
\title{Soft-constrained output feedback guaranteed cost equilibria in infinite-horizon uncertain linear-quadratic differential games}
\author[1]{Aniruddha Roy \thanks{roy9aniruddha@gmail.com}}
\author[1]{Puduru Viswanadha Reddy \thanks{vishwa@ee.iitm.ac.in}}
\affil[1]{Department of Electrical Engineering, \par Indian Institute of Technology--Madras, Chennai--600036, India} 
\date{ }
\maketitle

\begin{abstract} 
In this paper, we study infinite-horizon linear-quadratic uncertain differential games with an output feedback information structure. We assume linear time-invariant nominal dynamics influenced by deterministic external disturbances, and players' risk preferences are expressed by a soft-constrained quadratic cost criterion over an infinite horizon. We demonstrate that the conditions available in the literature for the existence of a soft-constrained output feedback Nash equilibrium (SCONE) are too stringent to satisfy, even in low-dimensional games. To address this issue, using ideas from suboptimal control, we introduce the concept of a soft-constrained output feedback guaranteed cost equilibrium (SCOGCE). At an SCOGCE, the players' worst-case costs are upper-bounded by a specified cost profile while maintaining an equilibrium property. We show that SCOGCE strategies form a larger class of equilibrium strategies; that is, whenever an SCONE exists, it is also an SCOGCE. We demonstrate that sufficient conditions for the existence of SCOGCE are related to the solvability of a set of coupled bi-linear matrix inequalities. Using semi-definite programming relaxations, we provide linear matrix inequality-based iterative algorithms for the synthesis of SCOGCE strategies. Finally, we illustrate the performance of SCOGCE controllers with numerical examples.
\end{abstract}
 
\section{Introduction}
\label{sec:Introduction}
Dynamic game theory provides a mathematical framework for modeling and analyzing multi-player decision-making situations that evolve over time. It has found successful applications in engineering, economics, management science, and biology, where these decision problems arise naturally; see \cite{isaacs:65,Engwerda:05,Basar:99,Dockner:00,Basar:18}. Specifically in engineering, dynamic game theory has been applied to address decision problems in, e.g., cyber-physical systems \cite{Zhu:15}, communication and networking \cite{Zazo:16,Basar:02}, autonomous vehicles \cite{Evens:22}, and smart grids  \cite{Chen:16,Maharjan:13,Mylvaganam:15}. In a dynamic game, the interaction environment of the players (or agents) is encapsulated by a state variable, and its evolution is described by a differential or difference equation. In a standard dynamic game, the dynamic model is assumed to be an exact representation of the interaction environment, and players use their strategies to achieve their objectives based on complete information about the state variable. However, real-world scenarios deviate from this standard framework in two significant ways. Firstly, obtaining an accurate dynamic model can be challenging, and optimization relying on highly detailed information from such a model may be impractical due to potential changes in the system's dynamics over time, which can be difficult to predict. Secondly, agents in real-world applications may not always have complete information about the state variable. For example, in distributed networked systems, agents can only access local state or output information related to their neighborhood, and they may not have direct access to the states but only to output information.

In the existing literature on dynamic games, one approach for addressing deviations of the first kind is through an uncertain dynamic game framework (see \cite{Van:2001thesis}, \cite{Van:03}, \cite{Engwerda:2006numerical}, \cite{Engwerda:2017robust}, \cite{Engwerda:22}). Here, a relatively simple model (a nominal model) is used to describe the dynamics. This approach is based on robust control theory \cite{Basar:08, Basar:03paradigms} and provides a way of capturing aversion to model risk in a dynamic game setting.
In this framework, it is assumed that the system dynamics is influenced by a deterministic disturbance, capturing discrepancies between the model and reality. Agents evaluate this disturbance according to their risk attitudes, which are reflected in their objective functions. They then manage this uncertainty by seeking strategies that optimize their objectives in a worst-case scenario.   In \cite{Van:2001thesis, Van:03}, two solution concepts are studied based on the norm of the disturbance signal under the feedback information structure:\footnote{The outcome of a differential game varies qualitatively based on the information available to the players during the decision-making process, also referred to as the information structure (see \cite{Basar:99}). In the open-loop information structure, the decisions of the players are functions of time and the initial state variable, while in the feedback information structure, they are functions of the current state variable. The latter information structure is   preferred as the decisions of the players are adapted to evolving state information.} the soft-constrained feedback Nash equilibrium, where the disturbance signal is unbounded but has finite energy, and the hard-bound constrained feedback Nash equilibrium, where the disturbance signal has a finite norm, implying that it is bounded. However, in all these works (including \cite{Engwerda:2006numerical}, \cite{Engwerda:2017robust}, \cite{Engwerda:22}), it is assumed that the agents have complete information about the state variable at all times during the course of the game. Consequently, these works do not address the second kind of deviation from the standard setup, which pertains to the situation where agents make decisions based on output information.

This paper is mainly concerned with the study of uncertain differential games with imperfect state or output feedback information structure. The contribution of our paper to the literature on dynamic games is twofold. Firstly, we consider a class of differential games with uncertain linear state dynamics and quadratic objectives defined over an infinite horizon. 
\textcolor{black}{We assume that system uncertainty is expressed by an additive disturbance signal to the nominal linear state dynamics, with external disturbances that are unbounded but have finite energy. Representative applications include smart grids, where agents manage energy distribution in the presence of short-term, unbounded fluctuations in renewable sources like wind and solar, which have finite energy over time. In financial markets, traders and automated agents respond to market shocks, such as price crashes, that are unbounded in the short term but exhibit finite volatility over longer periods. Similarly, in multi-robot systems, robots cooperate on tasks like exploration or mapping while facing unbounded environmental disturbances, such as impacts or terrain changes, with finite overall energy.}
When a state feedback information structure is assumed, \cite{Van:2001thesis, Van:03, Engwerda:22} show that the sufficient conditions for the existence of soft-constrained feedback Nash equilibria (SCFNE) are closely related to the existence of stabilizing solutions of a set of coupled algebraic Riccati equations (CARE). Following \cite{Engwerda:08}, with an output feedback information structure, the sufficient conditions for the existence of soft-constrained output feedback Nash equilibria (SCONE) require that the solutions of CARE satisfy additional structural constraints in addition to the stability requirements. We demonstrate that these sufficient conditions are too stringent to meet even for low-dimensional games. To address this problem, we introduce the notion of a soft-constrained output feedback guaranteed cost equilibrium (SCOGCE). This solution concept is inspired by the satisfaction equilibrium studied in \cite{Ross:06} in the context of static games. These strategies ensure that the worst-case individual costs of the players are upper-bounded by a given threshold (a design parameter) while retaining an equilibrium property. We demonstrate that SCOGCE strategies are a larger class of equilibrium strategies; that is, if SCONE strategies exist, they are contained in this class. Further, we derive sufficient conditions for the existence of SCOGCE strategies. The design of SCOGCE strategies uses techniques developed for suboptimal static output feedback controllers \cite{Leibfritz:01, Skelton:95, Skelton:17}. Our second contribution is related to the synthesis of SCOGCE strategies. In the existing literature, the SCFNE (or SCONE) strategies are obtained by first solving CARE, which involves solving \(Nn(n+1)/2\) multi-variable polynomial equations of degree 2 in \(Nn(n+1)/2\) variables, where \(N\) and \(n\) respectively denote the number of agents and the dimension of the state space. Though there exist iterative methods \cite{Li:1995, Mukaidani:2006_newton, Mukaidani:2004_lyapunov, Engwerda:07} for solving CARE, they are highly sensitive to initialization, and convergence to a stabilizing solution of CARE is not guaranteed. As a result, determining if a solution to CARE is stabilizing can only be verified ex-post, as demonstrated in \cite{Van:03, Engwerda:05, Engwerda:22,  Engwerda:2006numerical},  in other words, the required stability conditions are decoupled or not integrated into the existing iterative schemes for solving CARE. In this work, we demonstrate that the sufficient conditions for the existence of SCOGCE strategies result in a set of coupled matrix inequalities. Consequently, the synthesis of these strategies is based on checking the feasibility of these inequalities. We use semi-definite programming relaxation to transform these conditions into convex feasibility problems and develop a linear matrix inequality-based iterative method for synthesizing SCOGCE strategies. In our approach, the required stability conditions are incorporated within the iterative scheme, eliminating the need for ex-post verification.

This paper is organized as follows. In Section \ref{sec:preliminariesandproblemstatement}, we provide preliminaries and the problem statement. In Section \ref{sec:GC}, we derive some results pertaining to the suboptimal control problem of uncertain linear systems (single agent problem). Using these results, in Section \ref{sec:SCOGCE}, we introduce the notions guaranteed cost response and SCOGCE, and derive some properties of these equilibria. In Subsections \ref{sec:verification} and \ref{sec:synthesis}, we present results pertaining to the verification and synthesis of SCOGCE, respectively. In Subsection \ref{sec:algorithm}, we develop an algorithm for the computation of SCOGCE. In Section \ref{sec:SCSGCE}, we specialize these results to state feedback information structure. Section \ref{sec:numerical} illustrates the performance of SCOGCE-based controllers for the distributed control of networked multi-agent systems. Finally, conclusions are presented in Section \ref{sec:conclusions}.

\subsection{Novelty and differences with the existing literature} 
Infinite horizon uncertain linear quadratic differential games with deterministic disturbances were studied in \cite{Van:2001thesis,Van:03} using methods developed in robust control theory \cite{Basar:08,Basar:03paradigms}. Assuming a state feedback information structure, the authors propose soft-constrained and hard-bound constrained feedback Nash equilibria as solution concepts and provide sufficient conditions for their existence. For scalar uncertain games, \cite{Engwerda:2006numerical} proposes a numerical algorithm for computing soft-constrained Nash equilibria. In \cite{Engwerda:2017robust}, the author examines the existence of robust equilibria with an open-loop information structure. In weakly coupled large-scale systems, \cite{Mukaidani:2009soft} investigates the existence of soft-constrained equilibria for uncertain games with stochastic uncertainties. Recently, \cite{Engwerda:22} explores the design of min-max robust equilibrium controls in these uncertain differential games. In \cite{Jungers:08}, bounded type Nash controls were derived for linear systems with polytypic uncertainties. In all these studies, robust Nash equilibria are obtained under a state feedback information structure. Our work differs from these works by considering an output feedback information structure. 
The equilibrium concept presented in this paper is inspired by the satisfaction equilibrium studied in \cite{Ross:06} in static games. In the distributed control of multi-agent systems, agents use strategies based on local (or partial) state information. In this area, the works \cite{Mylvaganam:16, Cappello:18, Cappello:21distributed} study modeling frameworks that involve Nash and $\epsilon$-Nash equilibrium strategies. This paper differs from these works in both focus and solution concept, as we seek to obtain a broader class of equilibrium strategies that guarantee the individual worst-case costs of the agents remain below a given threshold.  This paper extends the work on uncertain differential games by  \cite{Van:2001thesis,Van:03,Engwerda:2006numerical} to a suboptimal setting. 
In our previous work \cite{Roy:22}, the guaranteed cost equilibrium was studied in networked linear quadratic differential games without external disturbances. This paper extends that work to a general class of uncertain differential games with an output feedback information structure.

\subsection{Notation} 
$\mathbb{R}^n$ ($\mathbb{R}^n_+$) denotes the set of $n\times 1$ real (positive) column vectors, $\mathbb{R}^{n\times m}$ denotes the set of $n\times m$ real matrices. $E^\prime$ denotes the transpose of a matrix or a vector $E$. ${I}_n$ denotes the $n\times n$ identity matrix. ${0}_{n\times m}$ denotes an $n\times m$ matrix with all its elements as zero. An $(n-1)$ tuple, $(1,2,\cdots,i-1,i+1,\cdots,n)$ is denoted by $-i$.
An $n$-tuple, $(S_1,S_2,\cdots,S_n)$ is also denoted by $(S_i,S_{-i})$ where  $S_{-i}=(S_1,\cdots,S_{i-1},S_{i+1},\cdots,S_n)$.  $\mathrm{col}\{e_1,e_2,\cdots,e_n\}$  denotes the single column vector or matrix obtained by column stacking the vectors or matrices $\{e_1,e_2,\cdots,e_n\}$. A positive definite (semi-definite) matrix $A$ is denoted by $A\succ 0$ ($A\succeq 0$), and $\sqrt{A}$ denotes its square root.  \textcolor{black} {For any two matrices $A, B \in \mathbb{R}^{n \times n}$, the relation $A \succ B$ ($A \succeq B$) denotes that the matrix $A - B$ is positive definite (positive semi-definite)}. The direct sum and Kronecker product of two matrices $A$ and $B$ are denoted by $A \oplus B$ and $A\bigotimes B$ respectively. A matrix $A$ is referred to as stable if the real part of all the eigenvalues of $A$ is negative. $||A||$ denotes a matrix norm of the matrix $A$. $S\bigtimes T$ denotes the Cartesian product of the sets $S$ and $T$.   $\partial S$ denotes the boundary of a set $S$. The space of $m$-dimensional real valued functions with finite energy (or quadratically integrable) on $[0,\infty)$ is denoted by $L_2^m([0,\infty)):=\left\{d(t)\in \mathbb R^m,~t\in[0,\infty)~|~\int_0^\infty d^\prime(t)d(t)dt <\infty \right\}$.   

\section{Preliminaries and problem formulation}
\label{sec:preliminariesandproblemstatement}
 
We denote the set of players  by $\N:=\{1,2,\cdots,N\}$. We assume that the dynamic interaction
environment of the players is affected by an external deterministic disturbance signal and is described by the following linear time-invariant dynamics 
\begin{subequations} \label{eq:system}
	\begin{align}
		\dot{x}(t)&=Ax(t)+\sum_{i\in \N} B_i u_i(t)+E d(t),~x(t)=x_0, \label{eq:sysdyn}\\ 
		y_i(t)&=C_ix(t),~i\in \N,\label{eq:sysobs}
	\end{align}%
where  $A\in \mathbb R^{n\times n}$,  $E\in \mathbb R^{n\times q}$, $B_i\in \mathbb R^{n\times m_i}$, and $C_i\in \mathbb R^{ s_i\times n}$ with $s_i \leq n$ and $\text{rank}(C_i)=s_i$.  We assume that the pairs $(A, B_i)$ and $(A, C_i)$ are stabilizable and detectable, respectively, for all $i \in \mathsf{N}$. Here, $x(t)\in \mathbb R^n$, $u_i(t)\in \mathbb R^{m_i}$, and $y_i(t)\in \mathbb R^{s_i}$ respectively denote the   state vector,  control action of Player $i$,  and output or observations of the state vector available to Player $i$  at time $t\in[0,\infty)$.  Further, $d(.)\in L_2^q([0,\infty))$ denotes $q$-dimensional deterministic disturbance signal, and  $x_0\in \mathbb R^n$ denotes the initial state vector.  Each Player $i\in \N$ chooses his controls so as to minimize the following soft-constrained quadratic cost functional defined over an infinite horizon
\begin{align}
	J_i(u_i,u_{-i},d):=\int_0^\infty \left(y_i^\prime(t) Q_i y_i(t)+u_i^\prime(t) R_i u_i(t)-d^\prime(t)D_id(t)\right) dt, \label{eq:playercost}
\end{align}
where   $Q_i\in \mathbb R^{s_i\times s_i}$, $Q_i\succeq 0$, $R_i\in \mathbb R^{m_i\times m_i}$, $R_i\succ0$ and $D_i\in \mathbb R^{q\times q}$, $D_i\succ 0$ for all $i\in \N$. We assume output feedback information structure, that is, players use static output feedback controls of the form 
\begin{align} 
	u_i(t)=F_i y_i(t),~F_i\in \mathbb R^{m_i\times s_i},~i\in \N \label{eq:fbcontrol}
\end{align} 
to achieve their objectives. Using this, the performance criterion \eqref{eq:playercost} is rewritten as
\begin{align}
	J_i(F_i,F_{-i},d)=\int_0^\infty \left(x^\prime(t) \left(C_i^\prime Q_i C_i +C_i^\prime F_i^\prime R_i F_i C_i\right) x(t)- d^\prime(t)D_i d(t) \right)dt. \label{eq:playerfbcost}
\end{align}
\textcolor{black}{To save on notation, we use the same letter $J_i$ 
to express both the original cost functional in \eqref{eq:playercost} and its reformulation  \eqref{eq:playerfbcost}}. Here, $x(t)$ evolves according to the closed-loop dynamics $\dot{x}(t)= \left(A+\sum_{i\in \N} B_i F_i C_i \right)x(t)+Ed(t)$, $x(0)=x_0$. In order to express players' desire for robustness to the unkonwn external disturbance, we assume that each player seeks to minimize his cost \eqref{eq:playerfbcost} with a worst-case disturbance, that is, Player $i\in \N$ aims to minimize the following modified
cost criterion
\begin{align}
	J_i^{\mathrm{sc}}(F_i,F_{-i}):=\sup_{d(.)\in L_2^q([0,\infty))} J_i(F_i,F_{-i},d). \label{eq:modifiedcost} 
\end{align}
\end{subequations} 
\begin{remark} In \eqref{eq:playerfbcost} the matrix $D_i$ represents Player $i$'s evaluation of the uncertainty. In particular, the term $d^\prime(t)D_id(t)$ appears with a negative sign. So, according to \eqref{eq:playerfbcost} and \eqref{eq:modifiedcost}, a large value of  $||D_i||$ indicates that Player  $i$'s worst-case disturbance signal will have a small magnitude. Consequently, Player $i$ does not anticipate significant deviations from the nominal dynamics in the environment, implying a risk-loving nature of the player. Conversely, a small value of $||D_i||$ indicates that Player $i$ is risk-averse.
\end{remark}  

\textcolor{black}{We assume that the dynamics \eqref{eq:sysdyn}-\eqref{eq:sysobs} and the objectives \eqref{eq:playerfbcost} are common knowledge among the players.} 
The equilibrium concept based on the adjusted costs \eqref{eq:modifiedcost} is given by the following
definition.
\begin{definition}[Soft-constrained output feedback Nash equilibrium]
	An $N$ tuple $(F^\star_i,F^\star_{-i}) $ is called a soft-constrained output feedback Nash equilibrium (SCONE) if for each Player $i \in \mathsf N$ the following inequality holds 
	\begin{align}
		& 	J_i^{\mathrm{sc}}(F^\star_i,F^\star_{-i}) \leq J_i^{\mathrm{sc}}(F_i,F^\star_{-i}), ~\forall F_i.
		\label{eq:SCNE}
	\end{align} 	
\end{definition}
For the full state feedback case, that is, when $s_i=n$ and $C_i=I_{n }$ for all $i\in \N$, \cite[Theorem 3.3]{Van:03} provides sufficient conditions for the existence of soft-constrained static state feedback Nash equilibrium. However, for the output feedback case, such conditions are not available. Nevertheless, following the study on output feedback Nash equilibria \cite[Theorem 3.2]{Engwerda:08}, an SCONE can be obtained using the sufficient conditions   \cite[Theorem 3.3]{Van:03}, and is related to the solvability of the following set of coupled algebraic Riccati equations 
\begin{align}
	\text{CARE}:~&\Bigl(A-\textcolor{black}{\sum_{j\in -i}}B_jR_j^{-1}B_j^\prime P_j\Bigr)^\prime P_i +P_i \Bigl(A-\textcolor{black}{\sum_{j\in -i}}B_jR_j^{-1}B_j^\prime P_j\Bigr)  \notag \\&\hspace{1in}+C_i^\prime Q_iC_i -P_iB_iR_i^{-1}B_i^\prime P_i + P_i ED_i^{-1}E^\prime P_i =0,~i\in \mathsf N. \label{eq:care1}
\end{align}
As the cost functions \eqref{eq:playerfbcost} are defined over an infinite horizon we require that players' strategies
should stabilize the closed-loop system. A stabilizing  solution of CARE is  an $N$-tuple $(P_i,P_{-i})$ of real symmetric matrices which satisfies \eqref{eq:care1} such that the following $N+1$ matrices 
\begin{align} 
A-\sum_{j=1}^{N} B_jR_j^{-1} B_j^\prime P_j,~ A- \sum_{j=1}^{N} B_j R_{j}^{-1} B_j^\prime P_j + ED_i^{-1}E^\prime P_i,~ i \in \mathsf{N}, \label{eq:stabcond}
\end{align} 
are stable. Following \cite[Remarks 2.2 and 2.3]{Engwerda:08}, if a  {stabilizing solution} of \eqref{eq:care1} can be found satisfying the following structural condition
\begin{align}
	B_i^\prime P_i (I-C_i^\prime (C_i C_i^\prime)^{-1}C_i)=0,~  i\in \mathsf N, \label{eq:structuralcond} 
\end{align} 
then the SCONE strategy of Player $i \in \mathsf N$ can be synthesized  as 
\begin{align} 
	F_i^\star=-R_i^{-1} B_i^\prime P_i C_i^\prime (C_i C_i^\prime)^{-1}. \label{eq:outfbstat}
\end{align}
The following example demonstrates that these sufficient conditions are quite restrictive, and fail to satisfy even for low-dimensional problems.
\begin{example}
	\label{ex:exp1}
	\begin{figure}[h] 
		\centering 
	\begin{tikzpicture}[scale=.25,>=latex', inner sep=1mm, font=\small]
		\tikzstyle{solid node}=[circle,auto=center,draw,minimum size=6pt,inner sep=2,fill=black!6]
		\tikzstyle{dedge} = [draw, blue!50, -,  line width=.5mm]
		\node (n1)[solid node] at (-5,0) {\scriptsize{$1$}};
		\node (n2)[solid node] at (0,0) {\scriptsize{$2$}};
		\node (n3)[solid node] at (5,0) {\scriptsize{$3$}}; 
		\path[dedge] (n1) --  (n2);
		\path[dedge] (n2) --  (n3);	 
	\end{tikzpicture}
		\caption{3-agent networked multi-agent system.}
\label{fig:3agentnetwork}
\end{figure}
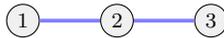
	Consider a 3-agent networked multi-agent system as shown in Figure \ref{fig:3agentnetwork}.
	The state dynamics associated with the agents is given by $\dot{x}_i(t)=x_i(t)+u_i(t)+ 0.1d(t)$, $x_i(t)\in \mathbb R$, $u_i(t)\in \mathbb R$, $i=1,2,3$, and $d\in L^1_2([0,\infty))$. The objectives of the players are given by
	$J_1=\int_0^\infty ( (x_1(t)-x_2(t))^2+u_1^2(t)-d^2(t))dt$, $J_2=\int_0^\infty ((x_2(t)-x_1(t))^2+(x_2(t)-x_3(t))^2 +u_2^2(t) -d^2(t) )dt$, and $J_3=\int_0^\infty ((x_3(t)-x_2(t))^2 +u_3^2(t)-d^2(t))dt$. The output variables of the players are given by $y_1(t)= \left[\begin{smallmatrix}x_1(t)&x_2(t) \end{smallmatrix}\right]^\prime$, $y_2(t)=\left[\begin{smallmatrix}x_1(t)&x_2(t)&x_3(t) \end{smallmatrix}\right]^\prime$ and $y_3(t)=\left[\begin{smallmatrix}x_2(t)&x_3(t) \end{smallmatrix}\right]^\prime$. The game parameters are
	$A=I_3$, $B_1=\left[\begin{smallmatrix} 1&0&0 \end{smallmatrix}\right]^\prime$,
	$B_2=\left[\begin{smallmatrix}0&1&0 \end{smallmatrix}\right]^\prime$,
	$B_3=\left[\begin{smallmatrix} 0&0&1\end{smallmatrix}\right]^\prime$,
	$E=0.1\left[\begin{smallmatrix}1&1&1\end{smallmatrix}\right]^\prime$,
    $Q_{1}=\left[\begin{smallmatrix} 1 &-1 \\ -1 &1 \end{smallmatrix}\right]$ , $Q_{2}=\left[\begin{smallmatrix} 1 &-1 &0 \\ -1 &2 &-1 \\ 0 &-1 &1 \end{smallmatrix}\right]$, 
	$Q_{3}=\left[\begin{smallmatrix} 1 &-1 \\ -1 &1 \end{smallmatrix}\right]$, $C_{1}=\left[\begin{smallmatrix} 1 &0 &0 \\0 &1 &0\end{smallmatrix}\right]$, $C_2=I_{3}$, $C_{3}=\left[\begin{smallmatrix}  0 &1 &0 \\ 0 &0 &1\end{smallmatrix}\right]$, 	$R_i=D_i=1$, for $i=1,2,3$. Solutions of CARE involve solving $18$ multi-variate polynomials of degree $2$ in $18$ variables. Upon imposing the  structural conditions \eqref{eq:structuralcond} the solutions of CARE must satisfy $P_1(1,3)=P_3(3,3)=0$. It can be verified that CARE does not admit solutions with these restrictions. So, it is not clear if  an SCONE exists for this game.
\end{example}  
\subsection{Problem statement}
\label{subsec:problem}
In this subsection, we outline the difficulties associated with computing an SCONE which serve as a  motivation for studying an alternate equilibrium concept.
\begin{enumerate}
\item  A solution for CARE requires solving $(Nn(n+1)/2)$ multivariate polynomial equations of degree $2$ in $(Nn(n+1)/2)$ variables. As mentioned in the introduction, although there exist iterative methods for solving CARE, they are highly sensitive to initial conditions, and convergence to a stabilizing solution is not guaranteed; refer to \cite{Engwerda:07}.

\item Once a solution to CARE is obtained, SCONE is determined by verifying if the solution satisfies the required $N+1$ stability conditions \eqref{eq:stabcond}. This implies that the stability requirements are not embedded but decoupled from the process of obtaining solutions to CARE; refer to \cite{Van:03,Engwerda:2006numerical,Engwerda:05}.

\item An additional complexity with SCONE is that the solution of a CARE must satisfy the structural condition \eqref{eq:structuralcond}. Similar to the stability requirement, this condition is decoupled, meaning it can only be verified after obtaining a solution to CARE; refer also to \cite{Engwerda:08}.
\end{enumerate}

Considering the limitations observed in the existence of SCONE, a natural question arises: Is there a broader class of output feedback strategies that satisfy an equilibrium property? Here, the term 'broader' refers to the requirement that SCONE strategies, if they exist, must be contained in this larger class of strategies. To address this question, we aim to develop an equilibrium concept based on a suboptimality criterion. In this framework, instead of players minimizing their worst-case costs, they seek output feedback strategies that guarantee their worst-case costs are lower than a given threshold.

\section{Suboptimal control of uncertain linear systems}
\label{sec:GC}
In this section, we consider a robust suboptimal linear quadratic control problem for  uncertain linear systems. First, we study the autonomous case and then synthesize   suboptimal feedback controls for uncertain linear systems with control inputs. Consider the autonomous system with external deterministic disturbance
\begin{subequations}\label{eq:suboptaut}
	\begin{align}
		\dot{x}(t) &= \bar{A} x(t) + Ed(t), ~ x(0)=x_{0}, \label{eq:dynamiceSingle}\\
		y(t) &=Cx(t),
	\end{align}
\end{subequations}
where  $\bar{A}\in \mathbb R^{n\times n}$ and $C\in \mathbb{R}^{s\times n}$ with $s\leq n$ and $\text{rank}(C)=s$. We assume that the pair $(A,C)$ is detectable. Here, $x(t)\in \mathbb R^n$, $y(t)\in \mathbb R^s$  and $d(t)\in \mathbb R^q$ denote the state, output and disturbance vectors respectively. We consider the following soft-constrained cost functional  
\begin{align}
	J (d):= \int_{0}^{\infty} \Bigl(y^\prime(t) \bar{Q} y(t) - d^\prime(t) D d(t)\Bigr)~dt,
	\label{eq:objectiveSingle}
\end{align}
where $\bar{Q}\in \mathbb R^{s\times s}$, $\bar{Q}\succeq 0$ and $D\in \mathbb R^{q\times q}$, $D\succ 0$. We seek to obtain conditions that ensure the worst-case cost, $\sup_{d \in L_2^{q}(0,\infty)} J(d)$, is lower than a given upper bound $\delta > 0$. The following lemma provides a sufficient condition.
\begin{lemma} \label{lem:gcaut}
	Consider the uncertain system \eqref{eq:suboptaut} with cost functional \eqref{eq:objectiveSingle}.    Let   $\delta >0$, and assume that there exist   $P\succ 0$ and $M\succ 0$ such that the following conditions hold
	\begin{subequations} 
		\begin{align}
			& \bar{A}^\prime P  + P \bar{A} + C^\prime \bar{Q} C+ PGP  \prec  0,
			\label{eq:LMI1disturbance}\\
			& x_{0}^{\prime}P x_{0} < \delta, \label{eq:LMI2upbound}\\
			& (\bar{A}+GP)^\prime M+M(\bar{A}+GP)\prec 0, \label{eq:LMI2stability} 
		\end{align}
		\label{eq:GCdisturbance}%
	\end{subequations} 
	where $G:=ED^{-1}E^\prime$. Then,  {$ \bar{A}$ and $(\bar{A}+GP)$ are stable},  and there exists a $\bar{d}(.) \in L_2^q([0,\infty))$ such that 
	$J( {d})\leq J(\bar{d})<\delta$ for all $d(.)\in L_2^q([0,\infty))$ where $\bar{d}(t)=D^{-1}E^\prime P e^{(\bar{A}+GP)t}x_0$.
\end{lemma}
\begin{proof} We write \eqref{eq:LMI1disturbance} as $\bar{A}^\prime P + P \bar{A} \prec -(C^\prime \bar{Q} C+ PGP)$. Then, since \textcolor{black}{$\bar{Q}\succeq 0$} and $D\succ 0$, and there exists a $P\succ 0$, it follows from the Lyapunov inequality that $\bar{A}$ is stable.
Next, as $d(.)\in L_2^q([0,\infty))$, it follows from Lemma \ref{lem:L2stability} (see Appendix) that $\lim_{t\rightarrow \infty} x(t) = 0$, which further implies that $\lim_{t\rightarrow \infty} y(t) = 0$. These observations imply that the cost \eqref{eq:objectiveSingle} is finite for all $d(.)\in L_2^q([0,\infty))$. Next, writing \eqref{eq:objectiveSingle} as 
	\begin{align}
		J(d)  
		&= \int_0^\infty \left(x^{\prime}(t) C^\prime\bar{Q}  C x(t) -d^\prime(t)Dd(t)\right) dt 
		- \int_{0}^{\infty} \frac{d}{dt}(x^{\prime}(t)P x(t))~dt  + \int_{0}^{\infty} \frac{d}{dt}(x^{\prime}(t)P x(t))~dt  .
		\label{eq:cost1}
	\end{align}
	As $ \lim_{t \rightarrow \infty} x(t) \rightarrow 0$,  we have $\int_{0}^{\infty} \frac{d}{dt}(x^{\prime}(t)P x(t))dt =\lim_{t \rightarrow \infty} x(t)^{\prime} P x(t) - x_0^{\prime}  P x_0=-x_0^\prime P x_0$. Using this in \eqref{eq:cost1}, and after completion of squares we get
	\begin{align}
		J(d)   
		&= x_0^{\prime}P x_0 + \int_{0}^{\infty} x^\prime(t)(\bar{A}^\prime P + P\bar{A}+ C^\prime \bar{Q}  C+ P G P)x(t) ~dt - \int_{0}^{\infty} ||d(t)-D^{-1}E^\prime P x(t)||^{2}_{D} ~dt. 
		\label{eq:cost2up}
	\end{align}
    From the finiteness of the cost for all $d\in L_2^q([0,\infty))$ and from the positive definiteness
	of $D$, we have that $J(d)$ is uniquely maximized at $\bar{d}(t)=D^{-1}E^\prime P\bar{x}(t) $, where $\bar{x}(.)$ is the state trajectory generated by $\dot{\bar{x}}(t) =(\bar{A}+GP) \bar{x}(t),~\bar{x}(0)=x_0$. From the Lyapunov inequality \eqref{eq:LMI2stability}, as there exists $M\succ0$, we have that $(\bar{A}+GP)$ is stable, which implies $\bar{d}(.)\in L_2^q([0,\infty))$. In particular, this implies $J(\bar{d})=\sup_{d \in L_2^{q}(0,\infty)} J(d)=\max_{d \in L_2^{q}(0,\infty)} J(d)$.
	
	Next, we denote by ${Q}^\dagger:= -(\bar{A}^\prime P + P\bar{A} + C^\prime\bar{Q} C + P G P)$,
	and write the worst-case cost as  
	\begin{align}
	 J(\bar{d}) = x_0^{\prime}P x_0 - \int_{0}^{\infty} \bar{x}^\prime(t){Q}^\dagger \bar{x}(t)~dt.
		\label{eq:cost3}
	\end{align}
Since  $\bar{A}+GP$ is stable and $Q^\dagger \succ 0$ (from \eqref{eq:LMI1disturbance}), we have that the quadratic functional $\int_{0}^{\infty} \bar{x}^\prime(t)Q^\dagger \bar{x}(t) dt$ is obtained uniquely 
	as $x_0^\prime Yx_0$ where $Y$ is a positive  definite solution of the Lyapunov equation
	$(\bar{A}+GP)^\prime Y+Y(\bar{A}+GP)=Q^\dagger$. Then, from \eqref{eq:LMI2upbound}, this implies $ 
	 J(\bar{d})=x_0^\prime Px_0 - x_0^\prime Y  x_0 \leq x_0^\prime P x_0 <\delta$.
\end{proof}
 
Next, we consider the controlled uncertain linear time-invariant system
\begin{subequations}\label{eq:uncertainlinsystem}
	\begin{align}
		\dot{x}(t)  &= A x(t) + Bu(t) + Ed(t), ~ x(0)=x_{0},
		\label{eq:dynamiceSingle1} \\
		y(t) &= Cx(t), 
	\end{align}
\end{subequations}
where   $B\in \mathbb R^{n\times m}$, and $u(t)\in \mathbb R^m$  denotes the control input. We assume that the pairs $(A,B)$ and $(A,C)$ are stabilizable and detectable respectively. We consider the following soft-constrained performance criterion
\begin{align}
	J(u,d,x_0):= \int_{0}^{\infty} \Bigl(y^\prime(t) Q y(t) + u^\prime(t)R u(t) - d^\prime(t) D d(t) \Bigr)~dt, 
	\label{eq:objectiveSingle1}
\end{align}
where  $R\in \mathbb R^{m\times m}$  with $R \succ 0$. 
We seek to obtain an output feedback control
$u(t)=Fy(t)$ which ensures that the worst-case cost $\text{sup}_{d\in L_2^q([0,\infty))}J(u,d,x_0)$ is lower than  $\delta >0$. The following lemma provides a sufficient condition for the existence 
of such a feedback control. 
\begin{lemma} Consider the uncertain controlled linear system \eqref{eq:uncertainlinsystem}  with the cost functional \eqref{eq:objectiveSingle1}.  Let $\delta >0$, and assume that there exist   $P\succ 0$, $M\succ 0$ and $F\in \mathbb R^{m\times s}$ such that the following conditions hold
	\begin{subequations} 
		\begin{align}
			& (A+BFC)^\prime P  + P (A+BFC) + C^\prime Q C+C^\prime F^\prime R FC + PGP  \prec 0, \label{eq:LMIdisturbance12}\\
			& x_{0}^{\prime}P x_{0} < \delta, \label{eq:LMIupbound2}\\
			& (A+BFC+GP)^\prime M+M(A+BFC+GP)\prec 0, \label{eq:LMIstability2}
		\end{align}
		\label{eq:GCdisturbance2}%
	\end{subequations} 
	where $G:=ED^{-1}E^\prime$.  {Then, $(A+BFC)$ and $(A+BFC+GP)$ are stable}, and  there exists a $\bar{d}(.)\in L_2^q([0,\infty))$ such that 
	$J(F, {d},x_0)\leq J(F,\bar{d},x_0)<\delta$ for all $d(.)\in L_2^q([0,\infty))$ 
 where \\ $\bar{d}(t)=D^{-1}E^\prime P e^{(A+BFC+GP)t}x_0$.
	\label{lem:GC_mainresult1} 
\end{lemma}
\begin{proof}
	The proof follows from Lemma \ref{lem:gcaut} by replacing $\bar{A}$ and $\bar{Q}$   with $A+BFC$ and $C^\prime Q C+ C^\prime F^\prime RF C$ respectively. 
\end{proof}
\begin{remark} Lemma \ref{lem:gcaut} extends the robust optimal control problem studied in \cite[Lemma 3.1]{Van:03}, for the autonomous case, to a suboptimal setting. Here, the additional stability condition \eqref{eq:LMI2stability} is crucial to ensure a finite upper bound for the worst-case cost, as demonstrated in \eqref{eq:cost3}. Similarly, Lemma \ref{lem:GC_mainresult1} provides a suboptimal counterpart to \cite[Lemma 3.2]{Van:03}.
\end{remark}
\section{Soft-constrained output feedback  guaranteed cost equilibrium}
\label{sec:SCOGCE}

In this section, we build upon the results obtained in the previous section (for the single-player case) and introduce the concept of a soft-constrained guaranteed cost function. Subsequently, we use this concept to introduce an equilibrium notion. Consider an $N$-tuple $(\delta_i, \delta_{-i}) \in \mathbb{R}_+^N$ representing a given cost profile, which is a design parameter. With the adjusted cost of Player $i$, as given by equation \eqref{eq:modifiedcost}, we define the set-valued worst-case guaranteed cost function for Player $i$ as follows:
\begin{align} 
	\mathsf g_i(F_{-i}):=\left\{F_i\in \mathbb R^{m_i\times s_i} ~|~ J^{\mathrm{sc}}_i(F_i,F_{-i}) <\delta_i\right\}.
	\label{eq:satisfactionfunction}
\end{align}
Here,   $\mathsf g_i(F_{-i})$ provides the set of all  output feedback strategies which ensure that Player $i$'s worst-case cost is below a given $\delta_i$, when other agents use their output feedback strategies \textcolor{black}{$F_{-i}\in \bigtimes_{j\in -i}\mathbb R^{m_j\times s_j} $}. Using this idea, we introduce an equilibrium concept in the next definition.
\begin{definition}[Soft-constrained output feedback guaranteed
	cost equilibrium] Let $ (\delta_i,\delta_{-i})\in \mathbb R^N_+$ be a given cost profile. 
	The strategy profile $  F^\circ: =(F^\circ_i,F^\circ_{-i}) $ is a \emph{soft-constrained output feedback guaranteed
		cost equilibrium} (SCOGCE) if for each Player $i\in \mathsf N$ the following condition holds true
	\begin{align}
		F_i^\circ\in \mathsf g_i(F_{-i}^\circ). \label{eq:fixedpoint}
	\end{align} 
	The set of all SCOGCE is given by
	\begin{align}
		\mathsf F^\circ:=\{(F_i,F_{-i}) \in   \bigtimes_{j\in \mathsf N} \mathbb R^{m_j\times s_j}~|~F_i\in \mathsf g_i(F_{-i}),~\forall i\in \mathsf N\}.
	\end{align}
	\label{def:NAGCE}
\end{definition}
\begin{remark}  We note that once the players are at an SCOGCE, no player has interest in deviating from the current strategy as each player has achieved a desired upper bound on his worst-case cost, in other words, players are reluctant to change their strategies once they are satisfied.
\end{remark} 
The next two results illustrate the properties of SCOGCE strategies. The first result is related to monotonicity property.
\begin{theorem}[Monotonicity]
	Let $(\delta_i,\delta_{-i})\in \mathbb R^N_+ $ and $(\bar{\delta}_i,\bar{\delta}_{-i})\in \mathbb R^N_+ $ be two cost profiles
	such that $\delta_i\leq \bar{\delta}_i$ for all $i\in \mathsf N$,
	and let $\mathsf F^\circ $ and $\bar{\mathsf F}^\circ $ denote the associated
	set of SCOGCE strategies respectively. Then,  $\mathsf F^\circ \subseteq \bar{\mathsf F}^\circ$. 
	\label{lem:LargeDelta}
\end{theorem}
\begin{proof} For any $(F_i^\circ,F_{-i}^\circ)\in \mathsf F^\circ$ we have that for every $i\in \N$, 
	$J^\mathrm{sc}_i(F_i^\circ,F_{-i}^\circ)<\delta_i$. Since $\delta_i\leq \bar{\delta}_i$, we have,
	for every $i\in \N$, $J^\mathrm{sc}_i(F_i^\circ,F_{-i}^\circ)<\bar{\delta}_i$, and this implies 
	$(F_i^\circ, F^\circ_{-i})\in \bar{\mathsf F}^\circ$.
\end{proof}
The next result shows that if an SCONE exists, then for a specific choice of  cost profile,  the SCONE is also an SCOGCE.
\begin{theorem}[Contains SCONE]  Let $(F_i^\star,F_{-i}^\star)$ be an SCONE for the uncertain LQDG. Then, $(F_i^\star,F_{-i}^\star)$ is also an SCOGCE considering the cost profile $\delta_i=J_i^\mathrm{sc}(F_i^\star,F_{-i}^\star)+\xi$,~$i\in \mathsf N$, for some $\xi>0$. 
	\label{thm:NashGCE}
\end{theorem}
\begin{proof} Since $\delta_i=J_i^\mathrm{sc}(F_i^\star,F^\star_{-i})+\xi$, we have $\mathsf g_i(F_{-i}^\star)\neq \emptyset$ as $J_i^\mathrm{sc}(F_i^\star,F^\star_{-i})<\delta_i$ for $i\in \N$.
	Further, from \eqref{eq:satisfactionfunction} and \eqref{eq:fixedpoint} this implies
	$F^\star_i\in \mathsf g_i(F_{-i}^\star)$ for every $i\in \N$. So, $(F_i^\star,F_{-i}^\star)$ is also an SCOGCE for the chosen cost profile.
\end{proof}
\begin{remark} Following Theorem \ref{thm:NashGCE}, we note that the set of guaranteed cost equilibrium strategies is larger; it contains the set of Nash equilibrium strategies.
\end{remark}
\begin{remark}
	The SCOGCE strategies can be interpreted as follows. Each Player $i\in \mathsf N$ privately reveals her guaranteed cost bound $\delta_i$ to a social planner or a principal. After knowing the cost profile $(\delta_i,\delta_{-i})$ of the players, the principal proposes a SCOGCE strategy profile  $(F_i^\circ,F_{-i}^\circ)$ to the players privately. Each player’s cost then achieves the desired upper bound by obeying the proposal when other players do the same. Thus, the SCOGCE is similar in spirit to a correlated equilibrium in static games \cite{Myerson:97}. 
\textcolor{black}{In this interpretation, it is sufficient for  player $i$ to know only her own threshold $\delta_i$, without needing knowledge of the thresholds of other players. Only the principal is required to have full knowledge of the thresholds for the synthesis of SCOGCE strategies. When the principal is not involved, each player is required to have knowledge of the thresholds of other players, in addition to her own, to synthesize her SCOGCE strategy. In this case, the cost profiles or thresholds of the players, as well as the notion of sub-optimality as rationality--where each player aims to achieve a cost below her threshold--should also be common knowledge among the players.}
\end{remark}
\subsection{Verification of SCOGCE} \label{sec:verification}
In this subsection, we   derive sufficient conditions to verify whether a given strategy profile is an SCOGCE for the uncertain differential game \eqref{eq:system}.
\begin{theorem}[Verification]   Let a cost profile $(\delta_i ,\delta_{-i} )\in \mathbb R^N_+$ and a strategy profile $F^\circ:=(F_i^\circ,F_{-i}^\circ)$ be given. Assume that for each $i\in \mathsf N$, there exist $P_i \succ 0$ and $M_i\succ 0$ such that the following conditions hold:
	\begin{subequations} 
		\begin{align} 
			& {A^\circ}^\prime P_i + P_i A^\circ + C_i^\prime Q_i C_i + C_i^\prime  F_i^{\circ \prime } R_i F_i^\circ C_i + P_i G_i P_i 
			\prec 0
			\label{eq:VR21}\\
			& x_{0}^{\prime}P_{i}x_{0} < \delta_i, \label{eq:VR22} \\
			& (A^\circ+G_iP_i)^\prime M_i+M_i(A^\circ+G_iP_i)\prec 0, \label{eq:VR23}
		\end{align}
		\label{eq:GCEoutput}%
	\end{subequations} 
	where $G_i:=E D_i^{-1}E^\prime$  and $A^\circ :=A+\sum_{i\in \mathsf N} B_i F^\circ_iC_i$. Then, $A^\circ$ and $A^\circ+G_iP_i$, $i\in \N$ are stable. Furthermore, $F^\circ$ is an SCOGCE, and the worst-case disturbance signal from Player $i$'s perspective is given by $\bar{d}_i(t)=D_i^{-1}E^\prime P_i e^{(A^\circ+G_iP_i)t}x_0$, $t\in[0,\infty)$, for $i\in \N$.
	\label{thm:GCEtestoutput} 
\end{theorem}
\begin{proof}
Firstly, we express \eqref{eq:VR21} as the Lyapunov inequality  ${A^\circ}^\prime P_i + P_i A^\circ \prec -(C_i^\prime Q_i C_i + C_i^\prime F_i^{\circ \prime } R_i F_i^\circ C_i + P_i G_i P_i)$. Given $Q_i\succeq 0$ and $R_i\succ 0$, along with the existence of $P_i\succ 0$ satisfying the Lyapunov inequality, we conclude that $A^\circ$ is stable. Additionally, the existence of $M_i\succ 0$ satisfying the Lyapunov inequality \eqref{eq:VR23} implies that $A^\circ + G_iP_i$ is stable. At $(F_i^\circ,F_{-i}^\circ)$, Player $i$'s cost is represented by:
\begin{align} 
	J_i(F_i^\circ,F_{-i}^\circ,d)=\int_0^\infty \left(x^\prime(t) \left(C_i^\prime Q_i C_i+C_i^\prime {F_i^\circ}^\prime R_i F_i^\circ C_i\right) x(t)-d^\prime(t) D_i d(t)\right) dt,
\end{align} 
where $x(t)$ satisfies $\dot{x}(t)=A^\circ x(t)+Ed(t)$ with $x(0)=x_0$. By replacing $\bar{A}$ and $\bar{Q}$ in Lemma \ref{lem:gcaut} with $A^\circ$ and $C_i^\prime Q_i C_i+C_i^\prime {F_i^\circ}^\prime R_i F_i^\circ C_i$ respectively, we obtain $J_i(F_i^\circ,F_{-i}^\circ,\bar{d})=J_i^\mathrm{sc}(F_i^\circ,F_{-i}^\circ)<\delta_i$, and
the worst-case disturbance from Player $i$'s perspective is given as $\bar{d}_i(t)=D_i^{-1}E^\prime P_i e^{(A^\circ+G_iP_i)t}x_0$ for $t\in[0,\infty)$. These conclusions hold for every player in $\N$, establishing that $(F_i^\circ, F_{-i}^\circ)$ is an SCOGCE.
\end{proof}
\begin{remark} We observe that the matrix inequality \eqref{eq:VR21} is quadratic in $P_i$.  Since $D_i \succ 0$,  using Schur complement (see \cite[Lemma 2.8]{Duan:13}),  \eqref{eq:VR21} is equivalent to the following LMI
	\begin{align}
		\begin{bmatrix}
			{A^\circ}^\prime P_i + P_i A^\circ + C_i^\prime Q_i C_i + C_i^\prime {F_i^{\circ}}^\prime R_i F_i^\circ C_i & P_i E \\
			E^\prime P_i & -D_i
		\end{bmatrix} \prec 0. \label{eq:VR21Schur}
	\end{align} 
Next, we note that \eqref{eq:VR23} is a bi-linear matrix inequality (BMI) in the variables $\{P_i, M_i, i \in \mathsf{N}\}$. As a consequence, feasibility of  \eqref{eq:VR23} cannot be verified using LMI solvers. However, using Lemma \ref{lem:Inq_upperbound}, we   obtain $
	(A^\circ + G_iP_i)^\prime M_i + M_i(A^\circ + G_iP_i) \preceq A^{\circ \prime}M_i + M_i A^\circ + \gamma_i P_i E D_i^{-1} E^\prime P_i + \frac{1}{\gamma_i} M_i E D_i^{-1} E^\prime M_i$, 
	for some $\gamma_i > 0$. This implies   the stricter condition $A^{\circ \prime} M_i + M_i A^\circ + \gamma_i P_i E D_i^{-1} E^\prime P_i + \frac{1}{\gamma_i} M_i E D_i^{-1} E^\prime M_i \prec 0$   guarantees the satisfaction of \eqref{eq:VR23}. Again, as $D_i \succ 0$, using  Schur's complement   the latter matrix inequality is equivalently written as the following LMI
	\begin{align}
		\begin{bmatrix}
			{A^\circ}^\prime M_i + M_i A^\circ & P_i E & M_i E \\
			E^\prime P_i & -\frac{1}{\gamma_i}D_i & 0_{q \times q} \\
			E^\prime M_i & 0_{q \times q} & -\gamma_i D_i
		\end{bmatrix} \prec 0. \label{eq:VR23Schur}
	\end{align} 
 {From the above, we can use the conditions \eqref{eq:VR21Schur}, \eqref{eq:VR22}, and the stricter condition \eqref{eq:VR23Schur} to verify if a strategy profile constitutes an SCOGCE   using   LMI solvers.}
\end{remark} 
\begin{remark}
Theorem \ref{thm:GCEtestoutput} extends our previous result \cite[Theorem 2]{Roy:22} in deterministic differential games with network constraints to uncertain differential games with output feedback information structure. In particular, Theorem \ref{thm:GCEtestoutput}  involves satisfaction of $3N$ matrix inequalities, whereas   \cite[Theorem 2]{Roy:22} involves $2N$ matrix inequalities. The additional $N$ inequalities \eqref{eq:VR23} are due to the stability requirement of the closed-loop system in the presence of disturbances.
\end{remark}

\subsection{Synthesis of SCOGCE: non-emptiness of the guaranteed cost response}
\label{sec:synthesis}
We observe that Theorem \ref{thm:GCEtestoutput} enables us to determine whether a given strategy profile is an SCOGCE for a specified cost profile. However, Theorem \ref{thm:GCEtestoutput} does not provide a method for synthesizing an SCOGCE. The synthesis of an SCOGCE relies on the non-emptiness of the guaranteed cost response \eqref{eq:satisfactionfunction}, and the following result provides the required sufficient condition.
\begin{theorem} Let $\delta_i>0$ and \textcolor{black}{$F_{-i}\in \bigtimes_{j\in -i}\mathbb R^{m_j\times s_j}$} be given.   Assume  there exist  $P_i\succ 0$, $M_i\succ 0$ and $F_i\in \mathbb R^{m_i\times s_i}$ such that 
	\begin{subequations}
		\begin{align}
			&(A_i+B_iF_iC_i)^\prime P_i + P_i (A_i+B_iF_iC_i) +C_i^\prime Q_iC_i + C_i^\prime F_i^\prime R_i F_iC_i +P_iG_iP_i \prec 0, \label{eq:LMI3}\\
			&x_0^\prime P_i x_0<\delta_i,\label{eq:LMI4}\\ 
			&(A_i+B_iF_iC_i+G_iP_i)^\prime M_i +M_i(A_i+B_iF_iC_i+G_iP_i)\prec 0, \label{eq:LMI5}
		\end{align}
		\label{eq:LMIlemoutput}%
	\end{subequations} 
	where $G_i:=E D_i^{-1}E^\prime$ and $A_i:=A+\textcolor{black}{\sum_{j\in -i}}B_j F_jC_j$. Then, the matrices $A_i+B_iF_iC_i$ and $A_i+B_iF_iC_i+G_iP_i$ are stable. Further, $\mathsf g_i(F_{-i})\neq \emptyset$, and 
	the worst-case disturbance signal from Player $i$'s perspective is given by $\bar{d}_i(t)=D_i^{-1}E^\prime P_i e^{(A_i+B_iF_iC_i+G_iP_i)t}x_0$ for $t\in[0,\infty)$. 
	\label{thm:GCoFi}  
\end{theorem} 
\begin{proof} The proof follows along the lines of the proof of Theorem \ref{thm:GCEtestoutput} and utilizes Lemma \ref{lem:GC_mainresult1}.
\end{proof}
We note that the sufficient condition \eqref{eq:LMI3} is a bi-linear matrix inequality (BMI) in the variables $(P_i,F_i)$, and \eqref{eq:LMI5} is a non-linear matrix inequality in the variables $(P_i,M_i,F_i)$, and their feasibility cannot be verified as convex feasibility problem; refer to \cite{Vanantwerp:2000}.  In the following we show that due to a special structure, the feasibility of \eqref{eq:LMI3} and \eqref{eq:LMI5} can be verified, using LMI solvers, with slightly stricter conditions, and consequently  provide an algorithm for the synthesis of $F_i\in \mathsf g_i (F_{-i})\neq \emptyset$. To this end,  we recall the following well-known lemma before stating the main result. 
\begin{lemma}(Projection theorem \cite[Theorem 2.3.12]{Skelton:17})
	Let $  S\in \mathbb R^{n\times m}$, $\text{rank}(  S)=m<n$,
	$T\in \mathbb R^{s\times n}$, $\text{rank}(T)=s<n$, and 
	$\Omega\in \mathbb R^{n\times n}$,~$\Omega=\Omega^\prime$, be given. Then, there exists $F\in \mathbb R^{m\times s}$
	satisfying $\Omega+S F T+ (S F T)^\prime  \prec 0$
	if and only if $\mathrm N_{S^\prime}^\prime \Omega \mathrm N_{S^\prime} \prec 0$ and $\mathrm N_T^\prime \Omega \mathrm N_T \prec 0$ hold, where $\mathrm N_{S^\prime}
	\in \mathbb R^{n\times (n-m)}$, $\mathrm N_T\in \mathbb R^{n\times (n-s)}$, denoting any matrices whose columns form orthonormal
	bases of the null spaces of $S^\prime$, $T$ respectively.
	\label{lem:Finsler}
\end{lemma} 
The next result provides a way to synthesize Player $i$'s feedback strategy $F_i\in \mathsf g_i(F_{-i})\neq \emptyset$ for a given 
$\delta_i>0$ and \textcolor{black}{$F_{-i} \in   \bigtimes_{j\in -i}\mathbb R^{m_j\times s_j}$}.
\begin{theorem}[Non-emptiness of the guaranteed cost response] \label{thm:LMIs} Let $\delta_i>0$ and \textcolor{black}{$F_{-i}\in \bigtimes_{j\in -i}\mathbb R^{m_j\times s_j}$} be given.    Consider the following sets
	\begin{subequations} 
		\begin{align}
			\mathsf X_i &:=\Big\{X\in \mathbb R^{n\times n} ~\big|~X\succ 0, ~  \mathrm N_{C^1_i}^\prime~ {\Omega}_i^1(X)~\mathrm N_{C_i^1}\prec 0 \Big\},\label{eq:XLMI2}\\
			\mathsf Y_i &:=\bigg\{Y \in \mathbb R^{n\times n}~\big|~Y\succ 0,~\begin{bmatrix}\delta_i &x_0^\prime\\x_0&Y\end{bmatrix}\succ  0,  ~ \mathrm N_{B_i^1}^\prime {\Omega}_i^2(Y)\mathrm N_{B_i^1}\prec 0\bigg\},\label{eq:YLMI2}\\ 
			\mathsf U_i(Z)&:=\Big\{U\in \mathbb R^{n\times n} ~\big|~U\succ 0,~\mathrm N_{C_i^2}^\prime \Omega_i^3(U;Z) \mathrm N_{C_i^2}\prec 0\Big\}, \label{eq:ULMI2}\\
			\mathsf V_i(Z)&:=\Big\{V\in \mathbb R^{n\times n} ~\big|~V\succ 0,~\mathrm N_{B_i^2}^\prime \Omega_i^4(V;Z) \mathrm N_{B_i^2}\prec 0\Big\}, \label{eq:VLMI2}\\
			{\Omega}^1_i(X)&:= \begin{bmatrix}A_i^\prime X+X A_i& C_i^\prime \sqrt{Q_i}& 0_{n \times m_{i}}& XE\\
				\sqrt{Q_i} C_i & -I_{s_{i}}& 0_{s_{i} \times m_{i}} & 0_{s_{i}\times q}\\ 0_{m_{i} \times n} & 0_{m_{i} \times s_{i}} & -R_i^{-1}&0_{m_{i} \times q}\\E^\prime X& 0_{q \times s_{i}} & 0_{q \times m_{i}}& -D_i \end{bmatrix}, \label{eq:OLMI12}\\ 
			{\Omega}^2_i(Y)&:=\begin{bmatrix}
				YA_i^\prime  +A_iY  & YC_i^\prime \sqrt{Q}_i & 0_{n \times m_{i}}&  E \\
				\sqrt{Q}_iC_iY & -I_{s_{i}} & 0_{s_{i} \times m_{i}} & 0_{s_{i} \times q}\\ 0_{m_{i} \times n} & 0_{m_{i} \times s_{i}}&-R_i^{-1}&0_{m_{i} \times q}\\ E^\prime & 0_{q \times s_{i}} &0_{q \times m_{i}}& -D_i \end{bmatrix}, \label{eq:OLMI22}  
			\\
			{\Omega}_i^3(U;Z)&:=\begin{bmatrix}A_i^\prime U+UA_i &ZE & UE\\
				E^\prime Z & - \frac{1}{\gamma_{i}} D_i & 0_{q \times q}\\
				E^\prime U  & 0_{q \times q} & -\gamma_{i} D_i \end{bmatrix}, \label{eq:OLMI32}\\
			{\Omega}_i^4(V;Z)&:=\begin{bmatrix}VA_i^\prime+A_iV &VZE &E\\
				E^\prime Z V & -\frac{1}{\gamma_{i}} D_i  & 0_{q \times q} \\
				E^\prime  & 0_{q \times q} & - \gamma_{i} D_i
			\end{bmatrix},\label{eq:OLMI42}
		\end{align}
		where $\gamma_i>0$ and $A_i =A+\textcolor{black}{\sum_{j\in -i}}B_j F_jC_j$. The matrices   
		$\mathrm N_{C_i^1}$,   $\mathrm N_{B_i^1}$, $\mathrm N_{C_i^2}$, and $\mathrm N_{B_i^2}$ denote matrices with orthonormal columns which span the null spaces
		of the matrices  $C_i^1:=\begin{bmatrix}C_i&0_{s_i \times s_i}&0_{s_i \times m_{i}}&0_{s_i \times q}\end{bmatrix}$, 
		$B_i^1:=\begin{bmatrix}B_i^\prime  &0_{m_i \times s_i}&I_{m_{i}}&0_{m_{i} \times q} \end{bmatrix}$,
		$C_i^2:=\begin{bmatrix}C_i&0_{s_i \times q}&0_{s_i \times q}\end{bmatrix}$, and
		$B_i^2:=\begin{bmatrix}B_i^\prime &0_{m_{i} \times q}&0_{m_{i} \times q} \end{bmatrix}$  respectively.  
		Define the sets
		\begin{align}
			\mathsf P_i&:=\Big\{P_i\in \mathbb R^{n\times n}  ~\Big|~P_i\succ 0,  ~P_i\in \mathsf X_i, ~P_i^{-1}\in \mathsf Y_i\Big\}, \label{eq:setP2_OF}\\
			\mathsf M_i(Z)&:=\Big\{M_i \in \mathbb R^{n\times n} ~\Big|~ M_i\succ 0,~M_i\in {\mathsf U}_i(Z),~M_i^{-1}\in \mathsf V_i(Z)\Big\}.\label{eq:setM2_OF}
		\end{align}
		When ${\mathsf P}_i\neq \emptyset$, we define the set 
		\begin{align}
			\mathsf {S}_i&:=\Big\{P_i\in \mathsf P_i~\Big|~ \mathsf{M}_i(P_i) \neq \emptyset\Big\}.\label{eq:setPs2_OF}
		\end{align} 
		If $\mathsf{S}_i\neq \emptyset$,  then  $\mathsf g_i(F_{-i})\neq \emptyset$. Further, for any $F_i\in \mathsf g_i(F_{-i})$, the matrices $A_i+ B_iF_i C_i$ and $A_i+B_iF_iC_i+G_iP_i$ are stable. Moreover, for any feasible $P_i \in \mathsf {S}_i$ and $M_i\in \mathsf {M}_i (P_i)$, an  $F_i\in \mathsf g_i(F_{-i})$ is obtained by solving the following LMIs
		\begin{align} 
			&\begin{bmatrix} (A_i+B_iF_iC_i)^\prime P_i+P_i(A_i+B_iF_iC_i)  +C_i^\prime Q_iC_i    & C_i^\prime F_i^\prime & P_i E \\
				F_iC_i & - R_i^{-1} &0_{m_{i} \times q}\\
				E^\prime P_i & 0_{q \times m_{i}}&- D_i \end{bmatrix} \prec 0,
			\label{eq:FGC2_OF}\\
			&(A_i+B_iF_iC_i+G_iP_i)^\prime M_i +M_i(A_i+B_iF_iC_i+G_iP_i)\prec 0.
			\label{eq:stability2_OF}  
		\end{align} 
		\label{eq:OFstab2}
	\end{subequations} 
	The worst-case disturbance signal from Player $i$'s perspective  is given by \\$\bar{d}_i(t)=D_i^{-1}E^\prime P_i e^{(A_i+B_iF_iC_i+G_iP_i)t}x_0$ for $t\in[0,\infty)$. 
	\label{thm:GCoi}
\end{theorem}
\begin{proof} As $Q_i\succeq 0$, writing $ {Q}_i=\sqrt{Q_i} (\sqrt{Q_i})^\prime$, using the Schur complement and the  notation provided in the theorem statement, the BMI  \eqref{eq:LMI3} is rewritten as  
	\begin{align}
		\Omega_i^1(P_i)+(\bar{B}_i^1)^\prime F_i C_i^1 + (C_i^1)^\prime F_i^\prime \bar{B}_i^1\prec 0,\label{eq:LMIint}
	\end{align} 
	where $\bar{B}_i^1= \begin{bmatrix} B_i^\prime P_i &0_{m_i\times s_i} &I_{m_i} &0_{m_i\times q}  \end{bmatrix}$.	
	From Lemma \ref{lem:Finsler}, the above BMI  \eqref{eq:LMIint} (in $(F_i,P_i))$ is feasible 
	if and only if the following LMIs (in $P_i$) 
	\begin{align}
		\mathrm N_{C_i^{1}}^\prime \Omega_i^{1}(P_i)  \mathrm N_{C_i^{1}} \prec 0, ~   {\mathrm {N}}_{\bar{B}_{i}^{1}}^\prime \Omega_i^{1}(P_i)   {\mathrm {N}}_{\bar{B}_{i}^{1}} \prec 0 
		\label{eq:projecction_lmi5}
	\end{align}
	are feasible. The nullspaces   $ {\mathrm{N}}_{\bar{B}_{i}^{1}}$ and  $\mathrm {N}_{B_{i}^{1}}$ are related as  
	$ {\mathrm {N}}_{\bar{B}_{i}^{1}} = (P_i^{-1}\oplus I_{s_{i}} \oplus I_{m_i}\oplus I_q )\mathrm {N}_{B_{i}^{1}}$, and using this, \eqref{eq:projecction_lmi5} (and as a result, \eqref{eq:LMI3}) is equivalently written as 
	\begin{align}
		\mathrm N_{C_i^{1}}^\prime \Omega_i^{1}(P_i)  \mathrm N_{C_i^{1}} \prec 0, ~  \mathrm {N}_{B_{i}^{1}}^\prime \Omega_{i}^2(P_i^{-1})  \mathrm {N}_{B_{i}^{1}} \prec 0.
		\label{eq:projection_lmi6}
	\end{align}
	Using Schur complement, the inequality \eqref{eq:LMI4} is written as
	\begin{align}
		\begin{bmatrix}
			\delta_i &x_0^\prime \\
			x_0 &P_i^{-1}
		\end{bmatrix} 
		\succ 0.
		\label{eq:LMI_ub2}
	\end{align}
	If $\mathsf{P}_i \neq \emptyset$, then from \eqref{eq:setP2_OF}, there is a $P_i \succ 0$ such that $P_i\in \mathsf{X}_i$ and $P_i^{-1}\in \mathsf{Y}_i$.  
	Using Lemma \ref{lem:Inq_upperbound}, we write  
		$(A_i+B_iF_iC_i+G_iP_i)^\prime M_i +M_i(A_i+B_iF_iC_i+G_iP_i) \preceq  
		(A_i+B_iF_iC_i)^\prime M_i + M_i(A_i+B_iF_iC_i) +
		\gamma_{i} P_i E D_i^{-1} E^\prime P_i   + \frac{1}{\gamma_{i}} M_i E D_i^{-1} E^\prime M_i$,
	for some $\gamma_i >0$. So, the satisfaction of the stricter inequality 
	$(A_i+B_iF_iC_i)^\prime M_i + M_i(A_i+B_iF_iC_i) +
		\gamma_{i} P_i E D_i^{-1} E^\prime P_i   + \frac{1}{\gamma_{i}} M_i E D_i^{-1} E^\prime M_i \prec 0$ 
	implies satisfaction of  \eqref{eq:LMI5}. 
	Using Schur complement, the latter matrix inequality equivalently written as 
	\begin{align}
		\begin{bmatrix}
			(A_i+B_iF_iC_i)^\prime M_i + M_i(A_i+B_iF_iC_i)  &P_i E &M_i E \\
			E^\prime P_i &-\frac{1}{\gamma_{i}} D_i &0_{q \times q} \\
			E^\prime M_i &0_{q \times q} & -\gamma_i D_i 
		\end{bmatrix}
		\prec 0.
		\label{eq:stab_oLMI3}
	\end{align}
	If the set $\mathsf{P}_i \neq \emptyset$, then for a feasible $P_i \in \mathsf{P}_i$,  \eqref{eq:stab_oLMI3} is a BMI in the variables  $(M_{i},F_{i})$, which is rewritten as 
	\begin{align}
		\Omega_i^{3}(M_i;P_i) + 
		( \bar{B}^2_i)^\prime
		F_i C^2_i  + (C^2_i)^\prime 
		F_i^\prime \bar{B}^2_i 
		\prec 0,
		\label{eq:proj_oLMI2}
	\end{align}
	where $\bar{B}_i^2 = \begin{bmatrix} B_i^\prime M_i &0_{m_i\times q} &0_{m_i\times q}  \end{bmatrix}$.
	Next,   using Lemma \ref{lem:Finsler}, the   BMI \eqref{eq:proj_oLMI2}  is feasible (in  $(F_i,M_i)$) if and only if the following LMIs   are feasible (in $M_i$)   
	\begin{align}
		\mathrm N_{C_{i}^2}^\prime  \Omega_i^{3}(M_i;P_i)  \mathrm N_{C_{i}^2} \prec 0, ~   {\mathrm {N}}_{\bar{B}_{i}^{2}}^\prime \Omega_i^3(M_i;P_i)  {\mathrm {N}}_{\bar{B}_{i}^{2}} \prec 0,
		\label{eq:proj_oLMI3}
	\end{align}
	where nullspaces   $ {\mathrm {N}}_{\bar{B}_{i}^{2}}$ and  $\mathrm {N}_{B_{i}^{2}}$ are related as  
	$	 {\mathrm {N}}_{\bar{B}_{i}^{2}} = (M_i^{-1}  \oplus I_{m_i}\oplus I_q )\mathrm {N}_{B_{i}^{2}}$. Using this, \eqref{eq:proj_oLMI3} (and as a result, \eqref{eq:proj_oLMI2}) is equivalently written as
	\begin{align}
		\mathrm N_{C_{i}^2}^\prime  \Omega_i^{3}(M_i;P_i)  \mathrm N_{C_{i}^2} \prec 0, ~\mathrm {N}_{B_{i}^{2}}^\prime \Omega_i^{4}(M_i^{-1};P_i)    \mathrm {N}_{B_{i}^{2}} \prec 0.
		\label{eq:proj_oLMI6}
	\end{align}
	If   $\mathsf{S}_i  \neq \emptyset$, then there is a $P_i \succ 0$ such that $P_i \in \mathsf{P}_i$ and $ \mathsf{M}_i(P_i) \neq \emptyset$, which further implies then there is a $M_i \succ 0$ such that $M_i \in \mathsf M_i (P_i) $. From \eqref{eq:setM2_OF}, this means, $M_i \succ 0$ satisfies $M_i\in \mathsf{U}_i (P_i)$ and  $M_i^{-1}\in \mathsf{V}_i(P_i)$. Then, following   Lemma \ref{lem:Finsler} there exists a $F_i \in \mathbb{R}^{m_{i} \times s_{i}}$ which is feasible for \eqref{eq:projection_lmi6} and \eqref{eq:proj_oLMI6} for any $P_i \in \mathsf{S}_i$ and $M_i \in \mathsf{M}_i(P_i)$. That is, there exist $P_i\succ0$, $M_i\succ 0$, and $F_i\in \mathbb R^{m_i\times s_i}$ which satisfy \eqref{eq:LMIlemoutput}. Then, from Theorem \ref{thm:GCoFi}, we have $\mathsf{g}_i(F_{-i}) \neq \emptyset$, and for any $F_i\in \mathsf g_i(F_{-i})$ we have that  $(A_i+ B_iF_iC_i) $ and $(A_i + B_iF_iC_i + G_i P_i) $ are stable. In particular, for a chosen $P_i \in \mathsf S_i$ and $M_i\in \mathsf M_i(P_i)$, the feedback strategy $F_i$ is synthesized as the one for which  LMIs \eqref{eq:FGC2_OF} and \eqref{eq:stability2_OF} are feasible.  
\end{proof}
\subsection{Algorithm for the synthesis of SCOGCE}
\label{sec:algorithm} 
  We note that Theorem \ref{thm:GCoi} provides a constructive approach that ensures Player $i$'s guaranteed cost response $\mathsf{g}_i(F_{-i})$ is non-empty for a given $(F_{-i},\delta_i)$. 
  In essence,  when the set $\mathsf{S}_i$, given by \eqref{eq:setPs2_OF}, is non-empty, then Player $i$'s output feedback strategy can be synthesized from \eqref{eq:FGC2_OF} and \eqref{eq:stability2_OF} such  that $J^\mathrm{sc}_i(F_i,F_{-i})<\delta_i$ and the matrices $A_i+B_iF_iC_i$ and $A_i+B_iF_iC_i+G_iP_i$ are stable. Using this, in this subsection, we propose an LMI-based iterative approach for computing an SCOGCE.

We note that the sets   $\mathsf{X}_i$ and $\mathsf{Y}_i$, defined using LMIs \eqref{eq:XLMI2} and \eqref{eq:YLMI2}, are open and convex. For computational purposes, we prefer closed sets, so we consider the following
closed ``$\epsilon$-approximation'' of the sets $\mathsf{X}_i$ and $\mathsf{Y}_i$ as
\begin{subequations} 
\begin{align} 
 	\mathsf X^\epsilon_i &:=\Big\{X\in \mathbb R^{n\times n} ~\big|~X\succ 0, ~  \mathrm N_{C^1_i}^\prime~ {\Omega}_i^1(X)~\mathrm N_{C_i^1}\preceq -\epsilon I_{n+s_i+m_i+q} \Big\},\label{eq:XLMI21}\\
 \mathsf Y^\epsilon_i &:=\bigg\{Y \in \mathbb R^{n\times n}~\big|~Y\succ 0,~\begin{bmatrix}\delta_i &x_0^\prime\\x_0&Y\end{bmatrix}\succeq  \epsilon I_{n+1},  ~ \mathrm N_{B_i^1}^\prime {\Omega}_i^2(Y)\mathrm N_{B_i^1}\preceq -\epsilon I_{n+s_i+m_i+q} \bigg\},\label{eq:YLMI21}
 \end{align}
where $\epsilon>0$. Consequently, the sets $\mathsf X_i^\epsilon$ and $\mathsf Y_i^\epsilon$  are closed, convex, 
and are contained respectively in $\mathsf X_i$ and $\mathsf Y_i$ (and also approach $\mathsf X_i$ and $\mathsf Y_i$ as $\epsilon\rightarrow 0$). We also approximate the convex, open sets $\mathsf U_i(Z)$ and $\mathsf V_i(Z)$ as follows 
 \begin{align} 
 \mathsf U^\epsilon_i(Z)&:=\Big\{U\in \mathbb R^{n\times n} ~\big|~U\succ 0,~\mathrm N_{C_i^2}^\prime \Omega_i^3(U;Z) \mathrm N_{C_i^2}\preceq -\epsilon I_{n+2q}\Big\}, \label{eq:ULMI21}\\
 \mathsf V^\epsilon_i(Z)&:=\Big\{V\in \mathbb R^{n\times n} ~\big|~V\succ 0,~\mathrm N_{B_i^2}^\prime \Omega_i^4(V;Z) \mathrm N_{B_i^2}\preceq -\epsilon I_{n+2q} \Big\}. \label{eq:VLMI21}
 \end{align} 
\end{subequations} 
Using the above approximated closed sets, the sets defined by \eqref{eq:setP2_OF}-\eqref{eq:setPs2_OF} are approximated as
\begin{subequations} 
		\begin{align}
		\mathsf P^\epsilon_i&:=\Big\{P_i\in \mathbb R^{n\times n}  ~\Big|~P_i\succ 0,  ~P_i\in \mathsf X^\epsilon_i, ~P_i^{-1}\in \mathsf Y^\epsilon_i\Big\}, \label{eq:setP2_OF1}\\
		\mathsf M^\epsilon_i(Z)&:=\Big\{M_i \in \mathbb R^{n\times n} ~\Big|~ M_i\succ 0,~M_i\in {\mathsf U}^\epsilon_i(Z),~M_i^{-1}\in \mathsf V^\epsilon_i(Z)\Big\}, \label{eq:setM2_OF1} \\
		\mathsf {S}^\epsilon_i&:=\Big\{P_i\in \mathsf P^\epsilon_i~\Big|~ \mathsf{M}^\epsilon_i(P_i) \neq \emptyset\Big\}.\label{eq:setPs2_OF1}
	\end{align} 
\end{subequations} 
Clearly, from the above $\mathsf S_i^\epsilon$ is a closed set, and it is contained in $\mathsf S_i$ and approaches $\mathsf S_i$ as $\epsilon \rightarrow 0$. Further, as $\mathsf S^\epsilon_i \subseteq \mathsf{P}^\epsilon_i$, to verify the feasibility of $\mathsf{S}^\epsilon_i \neq \emptyset$, we must verify the feasibility of $\mathsf{P}^\epsilon_i \neq \emptyset$. We note that, though the sets $\mathsf{X}^\epsilon_i$ and $\mathsf{Y}^\epsilon_i$, defined using LMIs \eqref{eq:XLMI21} and \eqref{eq:YLMI21}, are convex and closed, the set  $\mathsf P^\epsilon_i$, characterized by \eqref{eq:setP2_OF1}, is closed but not a convex set. As a result, the feasibility of $\mathsf P_i^\epsilon\neq \emptyset$ cannot be verified as a convex feasibility problem. 
The set $\mathsf P_i^\epsilon$ can be written equivalently using a bi-linear constraint as follows 
\begin{align}
	\mathsf {P}^\epsilon_i &:=\Bigl\{P_i, W_i \in \mathbb R^{n\times n}  ~\Big|~P_i\succ 0, W_i \succ 0,~ P_iW_i=I_{n},~P_i\in \mathsf{X}^\epsilon_i, ~W_i \in \mathsf {Y}_i^\epsilon\Bigr\}, \label{eq:setP2_bilinear}
\end{align}
Then, replacing the coupling constraint $P_iW_i=I_{n}$ using the semi-definite programming (SDP) relaxation
\cite{Leibfritz:01}
as $P_iW_i\succeq I_{n} \Leftrightarrow \left[\begin{smallmatrix}P_i& I_n \\I_n & W_i\end{smallmatrix}\right]\succeq 0$  provides a convex approximation of $\mathsf P_i^{\epsilon}$ as follows:  
\begin{align}
	\bar{\mathsf {P}}^\epsilon_i &:=\Bigl\{P_i, W_i \in \mathbb R^{n\times n}  ~\Big|~P_i\succ 0, W_i \succ 0,   \begin{bmatrix} P_i  &I_{n} \\ I_n  &W_i \end{bmatrix}
	\succeq 0,  ~ P_i\in \mathsf {X}^\epsilon_i, ~W_i \in \mathsf {Y}^\epsilon_i \Bigr\}. \label{eq:setP2_approx2}
\end{align}
In particular, we note that $(P_i,W_i)\in\partial \bar{\mathsf P}^\epsilon_i$ implies $P_iW_i=I_n$, and as a result  implies $\mathsf P^\epsilon_i=\partial \bar{\mathsf P}^\epsilon_i$. Once the feasibility of $\partial \bar{\mathsf P}^\epsilon_i\neq \emptyset$ is verified as a convex program, we need to verify
the feasibility of  $\mathsf M^\epsilon_i (P_i)\neq \emptyset$ for any $P_i\in \partial \bar{\mathsf P}^\epsilon_i$.
We note from \eqref{eq:ULMI2}, \eqref{eq:VLMI2}, and \eqref{eq:setM2_OF}, even though the sets
$\mathsf U^\epsilon_i(P_i)$ and $\mathsf V^\epsilon_i(P_i)$ are convex, the set $\mathsf M^\epsilon_i(P_i)$ is not a convex set. Again, using semi-definite relaxation we obtain a convex approximation of $\mathsf M^\epsilon_i(P_i)$ as
\begin{align}
	\bar{\mathsf M}^\epsilon_i(P_i) &:=\Bigl\{M_i, N_i \in \mathbb{R}^{n \times n} ~\Big| M_i\succ 0,N_i \succ 0, ~  \begin{bmatrix} M_i  &I_{n} \\ I_n  &N_i \end{bmatrix} \succeq 0,   ~M_i\in \mathsf {U}^\epsilon_i(P_i),~M_i^{-1}\in \mathsf {V}^\epsilon_i(P_i)\Bigr\},
	\label{eq:setMi_approx2}
\end{align}
where $\mathsf M^\epsilon_i(P_i)=\partial \bar{M}^\epsilon_i(P_i)$. Using the above, the set $\mathsf S_i$ is now written as 
\begin{align} \mathsf S^\epsilon_i=\left\{  P_i\in \partial \bar{\mathsf P}^\epsilon_i~|~ \partial \bar{\mathsf M}^\epsilon_i(P_i)\neq \emptyset   \right\}. \end{align}
So, the feasibility of $\mathsf S^\epsilon_i\neq \emptyset$
can be verified by first verifying $\partial \bar{\mathsf P}^\epsilon_i\neq \emptyset$ and 
then restricting the set  to the elements $P_i\in \partial \bar{\mathsf P}^\epsilon_i$ which ensure
the feasibility of $\partial \bar{\mathsf M}^\epsilon_i(P_i)\neq \emptyset $.  
To assess the feasibility of $\mathsf g_i(F_{-i}) \neq \emptyset$, for a given $(\delta_i,F_{-i})$, we introduce Algorithm \ref{alg:1}, named the Nested Sequential Linear Programming Matrix Method (NSLPMM). This method draws from the sequential linear programming matrix approach discussed in \cite[Algorithm 1]{Leibfritz:01}, originally devised for designing suboptimal static $\mathcal H_2/\mathcal H_\infty$ output feedback controllers.  The algorithm is structured into two stages, each targeting trace minimization problems. To begin, the initial stage involves solving the problem $\min \mathrm{tr}(P_iW_i)$ where $(P_i,W_i)\in \bar{\mathsf P}^\epsilon_i$. If this yields a boundary solution, denoted as $P_i\in \partial \bar{\mathsf P}^\epsilon_i$, the next step addresses the second problem, $\min \mathrm{tr}(M_iN_i)$, with $(M_i,N_i)\in \bar{\mathsf M}^\epsilon_i(P_i)$.
It's crucial to emphasize that both these problems are well-posed \cite{Leibfritz:01}, as their objectives are bounded below by $n$. The first stage, spanning Steps 1-5 and 31-33, adopts an iterative strategy to solve the initial minimization problem. The process starts by determining $(P^0_i,W^0_i) \in \bar{\mathsf P}^\epsilon_i$, which is an LMI feasibility problem. In each iteration ($k=0,1,\cdots$), an SDP minimization problem, featuring a linear objective and LMI constraints, is solved; see Step 4. 
Following \cite[Theorem 3.9]{Leibfritz:01}, this iterative process ensures that $(P_i^k, W_i^k)$ lies within $\bar{\mathsf P}^\epsilon_i$ for all $k=0,1,2,\cdots$ and that the sequence $(P_i^k, W_i^k)$ converges within $\bar{\mathsf P}^\epsilon_i$ towards a minimum of the problem $\min \mathrm{tr}(P_iW_i)$. Step 5 verifies if this minimum is achieved. If the solution $(P_i^k, W_i^k)$ attains the lower bound where $\mathrm{tr}(P_i^kW_i^k)=n$, reached at the boundary $(P_i^k, W_i^k)\in \partial \bar{\mathsf P}^\epsilon_i$, then $\mathsf P^\epsilon_i \neq \emptyset$ and the algorithm proceeds to the second stage. Step 6 checks if the minimum indeed reaches the lower bound. If not, the algorithm stops, and the guaranteed cost response $\mathsf g_i(F_{-i})$ remains empty. 
In the second stage, beginning with the solution $P_i \in \mathsf P^\epsilon_i$ obtained in the first stage, the approach to solving the second minimization problem follows a similar process as described in Steps 8-11 and 22-24. If the minimum does achieve the lower bound, specifically when $\mathrm{tr}(M_i^lN_i^l)=n$ and $(M_i^l, N_i^l)\in \partial \bar{\mathsf M}_i^\epsilon(P_i)$, then $\mathsf S_i^\epsilon \neq \emptyset$, and by Theorem \ref{thm:LMIs}, this implies that the guaranteed cost response $\mathsf g_i(F_{-i})$ is non-empty. Then, Player $i$'s feedback strategy, ensuring the cost $J_i^{\mathrm{sc}}(F_i,F_{-i}) < \delta_i$, is synthesized from \eqref{eq:FGC2_OF}-\eqref{eq:stability2_OF}. However, if the minimum fails to attain the lower bound, the algorithm terminates, and the guaranteed cost response $\mathsf g_i(F_{-i})$ remains empty.

Next, we propose an iterative algorithm to find an SCOGCE which satisfies the sufficient conditions \eqref{eq:GCEoutput} stated in  \textcolor{black}{Theorem \ref{thm:GCEtestoutput}}. Algorithm \ref{alg:2}, referred to as sequential guaranteed cost response, is based on sequential best response algorithm which is widely used in the game theory literature \cite{Engwerda:2006numerical}, \cite{Tembine:2011book}; see also \cite{Roy:22}, \cite{Veetaseveera:19}. The algorithm is initialized with $(F_i, F_{-i})$, ensuring the matrix $(A + \sum_{i \in \mathsf{N}}B_iF_iC_i)$ remains stable,
and this is achieved by solving the   feasibility of Lyapunov inequality $\left(A+ \sum_{i=1}^{N}B_i F_i C_i \right)^\prime P + P \left(A+ \sum_{i=1}^{N} B_i F_i C_i\right) \prec 0$ in the variables   $(F_i,F_{-i})$ for a given $P\succ 0$. In Step 2, Algorithm \ref{alg:1} is used to verify $\mathsf S_i^\epsilon \neq \emptyset$ for a fixed $F_{-i}$. In Step 3, $F_i$ is obtained using \eqref{eq:FGC2_OF}-\eqref{eq:stability2_OF}  for a feasible $P_i\in \mathsf S_i^\epsilon \neq \emptyset$. In Step 10, the player index is updated to the next player and the procedure is repeated until an SCOGCE is found. The algorithm does not stop immediately even if Player $i$'s strategy is not found in Step 3. Only when all the players consecutively fail to obtain the strategy in Step 3 the algorithm is halted. It is possible that the algorithm can cycle forever without stopping involving updates of a few players' strategies. It is well known \cite{Tembine:2011book} that the best response algorithm converges only for a certain classes of games. In our work, an added complexity is due to the non-convex nature of the problem (in Theorem \ref{thm:GCoi}) requiring NSLPMM algorithm to verify $\mathsf S_i^\epsilon \neq \emptyset$. This shortcoming is inherent in the design of static output feedback suboptimal controllers; see \cite{Leibfritz:01} and \cite{Iwasaki:94}. For this reason, we do not have theoretical guarantees for the convergence of Algorithm \ref{alg:2}.
\begin{remark} As mentioned in the problem formulation, the iterative methods for solving the CARE in \cite{Van:03, Engwerda:2006numerical, Engwerda:22} do not guarantee the satisfaction of the necessary stability conditions,  that is, the stability conditions have to be verified separately after a solution to CARE is obtained. In the proposed approach, the required stability conditions are integrated at every step in the iterative process (refer to \eqref{eq:LMI3} and \eqref{eq:LMI5}).
\end{remark} 
\begin{remark}
	In \cite{Roy:22}, the guaranteed cost equilibrium concept was studied in distributed control for networked heterogeneous multi-agent systems without external disturbances. This paper extends this concept to a general class of uncertain linear quadratic differential games, and as a result can also be used to design robust distributed control of networked heterogeneous multi-agent systems. 
\end{remark}
\begin{remark} We note that the design of SCOGCE controllers is centralized but their implementation is distributed. The synthesis of SCOGCE can be made independent of global initial state $x_0$ by assuming that
	$x_0\in  \{z\in \mathbb R^n~|~z'z<\alpha\}$, where $\alpha>0$ is a scaling factor.
	Then, the conditions \eqref{eq:VR22} and \eqref{eq:LMI5} are replaced with $P_i\prec \tfrac{1}{\alpha} \delta_i I_n$, and   \eqref{eq:setP2_OF} is replaced with 
	$	\mathsf P_i:=\{P_i\in \mathbb R^{n\times n} ~|~0\prec P_i\prec \tfrac{1}{\alpha} \delta_i I_n,  ~P_i\in \mathsf X_i, ~P_i^{-1}\in \mathsf Y_i\}$.
\end{remark}
\begin{remark}
	The selection of the cost profile $(\delta_i,\delta_{-i})\in \mathbb R^N_+$ plays a central role in the existence of SCOGCE. The sufficient conditions stated in Theorems  \ref{thm:GCoFi} and \ref{thm:GCEtestoutput} are highly interrelated; that is, non-emptiness of guaranteed cost response of Player $i$ for a cost estimate $\delta_i$ also depends on $\delta_{-i}$, the cost estimates of other players. Following the monotonicity property of SCOGCE (see Theorem \ref{lem:LargeDelta}), one approach to get these estimates would be to use large values of the cost profile, in Algorithms \ref{alg:1} and \ref{alg:2}, and slowly decrease them until an SCOGCE cannot be found.
\end{remark}
\begin{remark}  In the game theory literature, it is well known that an equilibrium outcome is inefficient, and players
	can do better in lowering their costs through cooperation. The Price of Stability (PoS) is a measure used to quantify the degradation in efficiency caused by the game-theoretic constraint of stability associated with an equilibrium; see  \cite{Nisan:07, Anshelevich:08} in static games and \cite{Zhu11} in differential games. It is defined as the ratio of the total game cost in equilibrium to that in cooperation. The PoS associated with an SCOGCE  strategy profile $F^\circ:=(F_i^\circ,F_{-i}^\circ)\in \mathsf F^\circ$ is given by 
	\begin{align}
		\textcolor{black}{\mathrm{PoS}(F^\circ)}=\frac{\sum_{i=1}^n J^\mathrm{sc}_i(F_i^\circ,F_{-i}^\circ)}{J^\text{Co}},
		\label{eq:PoSdef}
	\end{align}
	where $J^\text{Co}$ is the total game cost in cooperation, and is obtained as the optimal cost associated with the robust optimal control problem
	\begin{align} 
		J^\text{Co}:=\min_{(F_i,F_{-i})} \sup_{d\in L_2([0,\infty))} \sum_{i=1}^N J_{i}(F_i,F_{-i},d),
		\label{eq:tgcostcoop}
	\end{align} 
	where the state and the output variables evolve  according to \eqref{eq:system}. 	
	\textcolor{black} {Whenever there are multiple SCOGCEs, the PoS, given by \eqref{eq:PoSdef}, can be used to rank these equilibria within the set $\mathsf{F}^\circ$. In particular, the equilibrium with the minimum PoS is preferred, as its efficiency is closest to the cooperative outcome; see also Section \ref{sec:numerical}. However, multiple equilibria may yield the same minimum PoS, and in such cases, we currently do not have a method to further refine them.} 
	\label{rem:PoS}
\end{remark} 
\vskip1ex 
\begin{algorithm}[H]
	\captionsetup{labelfont={sc,bf}, labelsep=newline}
	\setstretch{1}
	\caption{SCOGCE: Non-emptiness of the guaranteed cost response using NSLPMM} 	\label{alg:1}
	Determine initial $(P_{i}^{0}, W_{i}^{0}) \in \bar{\mathsf{P}}^\epsilon_{i}$;
	\\
	\For{$k=0, 1, 2, \cdots $}{
		Determine $(S^{k}_1,T^{k}_1)$ as the unique solution of  \\
		$\min \mathrm{tr}(P_{i}W_{i}^{k}+P_{i}^{k}W_{i})$ \text{s.t.} $(P_{i},W_{i}) \in \bar{\mathsf{P}}^\epsilon_{i}$ 
		\\
		\eIf{ $\mathrm{tr}(S^{k}_1W_{i}^{k}+P_{i}^{k}T^{k}_1) = 2 \mathrm{tr}(P_{i}^{k}W_{i}^{k})$
		}	
		{ \eIf{$\mathrm{tr}(P_k^iW_i^k)=n$} {
				Found $P_i^k \in \mathsf P_i^\epsilon= \partial \bar{\mathsf P}^\epsilon_i$. Set $P_i=P_i^k$\\
				Determine initial $(M_i^0,N_i^0)\in \bar{\mathsf M}^\epsilon_i(P_i)$\\
				\For{$l=0,1,2\cdots $}{
					Determine $(S^l_2,T^l_2)$ as the unique solution of \\
					$\min \mathrm{tr}(M_iN_i^l+M_i^lN_i)$ s.t. $(M_i,N_i)\in \bar{\mathsf M}^\epsilon_i(P_i)$\\
					\eIf{$\mathrm{tr}(S^{l}_2N_{i}^{l}+M_{i}^{l}T^{l}_2) = 2 \mathrm{tr}(M_{i}^{l}N_{i}^{l})$}{
						\eIf{$\mathrm{tr}{(M_i^lN_i^l)}=n$}{
							Found $M_i^l \in \mathsf M^\epsilon_i(P_i)= \partial \bar{\mathsf M}^\epsilon_i(P_i)$\\
							Set $M_i=M_i^l$\\ Found $(P_i,M_i)$ s.t. $P_i\in \partial \bar{\mathsf{P}}^\epsilon_i$ and $M_i\in \partial \bar{\mathsf M}^\epsilon_i(P_i)$, this implies 
							$\mathsf g_i(F_{-i}) \neq \emptyset$\\
							Stop. }{Stop. $\mathsf g_i(F_{-i})=\emptyset$}}{
						Compute $\beta \in [0,1]$ by solving \\
						$\displaystyle \min_{\beta  \in [0,1] } \mathrm{tr} ((M_{i}^l+\beta(S^l_2-M_{i}^l))(N_{i}^l+\beta (T^l_2-N_{i}^l))) $ 
						\\
						Set $M_i^{l+1}=(1-\beta )M_i^l+\beta S^l_2$, $N_i^{l+1}=(1-\beta)N_{i}^l+\beta T^l_2$}
				} 
			} 
			{Stop. $\mathsf g_i(F_{-i})=\emptyset$}}
		{Compute $\alpha \in [0,1]$ by solving \\
			$\displaystyle \min_{\alpha \in [0,1] } \mathrm{tr} ((P_{i}^{k}+\alpha (S^{k}_1-P_{i}^{k}))(W_{i}^{k}+\alpha  (T^{k}_1-W_{i}^{k}))) $ 
			\\
			Set $P_i^{k+1}=(1-\alpha )P_i^{k}+\alpha S^{k}_1$, $W_i^{k+1}=(1-\alpha )W_{i}^{k}+\alpha  T^{k}_1$}  
	} 
\end{algorithm}
{\begin{algorithm} [H]
		 \setstretch{1}
		\caption{SCOGCE: Sequential guaranteed cost response}
		\label{alg:2}
		\DontPrintSemicolon  
		\KwData{ given   $(\delta_i,\delta_{-i})\in \mathbb R_+^N$ and initial stabilizing guess $(F_i,F_{-i})$;} 
		\KwIn{$i=1$; \tcp*[l]{\small{player index}}}
		\KwIn{$j=0$; \tcp*[l]{\small{number of players who fail to satisfy Theorem \ref{thm:GCoi} consecutively}}}
		\While{$(F_i,F_{-i})$ is not an SCOGCE (verified using Theorem \ref{thm:GCEtestoutput})}
		{
			\eIf {$ \mathsf{S}^\epsilon_i  \neq \emptyset$ and $ \mathsf{M}^\epsilon_i (P_i) \neq \emptyset$ exist satisfying Theorem  \ref{thm:GCoi} given $F_{-i}$, using Algorithm \ref{alg:1}}
			{ Update $F_i$ using \eqref{eq:FGC2_OF} and \eqref{eq:stability2_OF} for a feasible $P_{i} \in \mathsf{S}_{i}^{\epsilon}$ and $M_i \in \mathsf{M}^\epsilon_i (P_i)$; \;
				$j=0$;
			}
			{$j=j+1$;\; 
				\If{$j \geq N$}{\KwResult{Stop, cannot find SCOGCE.}}}
			$i=(i~\text{mod}~N)+1$ \tcp*[l]{\small{update player index}} }
		\KwResult{Found an SCOGCE}
\end{algorithm}}
\section{Soft-constrained state feedback guaranteed cost equilibrium}
\label{sec:SCSGCE}
 In this section, we specialize the results obtained in the previous section to the situation where all the players have access to complete state information, that is, when $C_i=I_n$ and $y_i(t)=x(t)$ for all $i \in \mathsf N$ and $t \in [0,\infty)$. The SCOGCE in this case is referred to as a soft-constrained state feedback guaranteed cost equilibrium (SCSGCE), and we seek to obtain sufficient conditions for the existence of an SCSGCE. To this end, firstly, the verification of an SCSGCE can be obtained by replacing \eqref{eq:VR21} in Theorem \ref{thm:GCEtestoutput} with ${A^\circ}^\prime P_i + P_i A^\circ + Q_i +  F_i^{\circ \prime } R_i F_i^\circ + P_i G_i P_i 
\prec 0$. Secondly, to synthesize SCSGCE strategies, we need to ensure that the guaranteed cost response of each player is non-empty. By specializing the matrix inequalities specializing the matrix inequalities \eqref{eq:LMIlemoutput} in Theorem \ref{thm:GCoFi}, we obtain the required sufficient conditions as follows:
\begin{subequations}
\begin{align} &(A_i+B_iF_i)^\prime P_i + P_i (A_i+B_iF_i) + Q_i+ F_i^\prime R_i F_i+P_iG_iP_i \prec 0, \label{eq:SCSGCE1}\\
	& x_0^\prime P_i x_0<\delta_i \label{eq:SCSGCE2},\\
	&(A_i+B_iF_i+G_iP_i)^\prime M_i +M_i(A_i+B_iF_i+G_iP_i)\prec 0.\label{eq:SCSGCE3}
	\end{align} 
\end{subequations} 
In the next corollary, we specialize Theorem \ref{thm:GCoi} to synthesize Player $i$'s feedback strategy $F_i\in \mathsf g_i(F_{-i})\neq \emptyset$ for a given $\delta_i>0$ and \textcolor{black}{$F_{-i} \in \bigtimes_{j\in -i}\mathbb R^{m_j\times n} $}. 
\begin{corollary} [Non-emptiness of the guaranteed cost response]   Let $\delta_i>0$ and 
	\textcolor{black}{$F_{-i}\in \bigtimes_{j\in -i}\mathbb R^{m_j\times n}$} be given.    Consider the following sets  
	\begin{subequations} 
		\begin{align}
			\mathsf Y_i &:=\Big\{Y\in \mathbb R^{n\times n} ~\big|~   Y\succ 0, ~\begin{bmatrix}\delta_i &x_0^\prime\\x_0& Y\end{bmatrix}\succ  0,~  \mathrm N_{B_i^3}^\prime {\Phi}^1_i(Y) \mathrm N_{B_i^3}\prec 0 \Big\},\label{eq:XLMI2_SF}\\
			\mathsf V_i(Z) &:=\bigg\{V \in \mathbb R^{n\times n}~\big|~ V\succ 0, \mathrm N_{B_i^4}^\prime{\Phi}^2_i(V;Z) \mathrm N_{B_i^4}\prec 0\bigg\},\label{eq:YLMI2_SF} \\ 
			{\Phi}^1_i(Y)&:= \begin{bmatrix}YA_i^\prime+ A_iY+G_i& Y \sqrt{Q_i}& 0_{n \times m_{i}}\\
				\sqrt{Q_i} Y& -I_{n} & 0_{n \times m_{i}}  \\ 0_{m_{i} \times n} & 0_{m_{i} \times n}  & -R_i^{-1}  \end{bmatrix}, \label{eq:OLMI12_SF}\\ 
			{\Phi}^2_i(V;Z)&:=\begin{bmatrix}VA_i^\prime+A_iV+ \tfrac{1}{\gamma_i} G_i &VZE \\
				E^\prime Z V & -\frac{1}{\gamma_{i}} D_i\end{bmatrix},
		\end{align}
		where $\gamma_i>0$ and $A_i =A+\textcolor{black}{\sum_{j\in -i}}B_j F_j$. The matrices   
		$\mathrm N_{B_i^3} $  and $\mathrm N_{B_i^4} $ denote matrices with orthonormal columns which span the null spaces
		of the matrices  
		$B_i^3:=\begin{bmatrix}B_i^\prime  &0_{m_{i} \times n}&I_{m_{i}}\end{bmatrix}$,
		and
		$B_i^4:=\begin{bmatrix}B_i^\prime &0_{m_{i} \times q}\end{bmatrix}$  respectively.  
		Define the sets
		\begin{align}
			\mathsf P_i&:=\Big\{P_i\in \mathbb R^{n\times n}  ~\Big|~P_i\succ 0,  ~P_i^{-1}\in \mathsf Y_i\Big\}, \label{eq:setP2}\\
			\mathsf M_i(Z)&:=\Big\{M_i \in \mathbb R^{n\times n} ~\Big|~ M_i\succ 0,~M_i^{-1}\in \mathsf V_i(Z)\Big\}.\label{eq:setM2}
		\end{align}
		When ${\mathsf P}_i\neq \emptyset$, we define the set 
		\begin{align}
			\mathsf {S}_i&:=\Big\{P_i\in \mathsf P_i~\Big|~ \mathsf{M}_i(P_i) \neq \emptyset\Big\}.\label{eq:setPs2} 
		\end{align} 
		If $\mathsf{S}_i\neq \emptyset$,  then  $\mathsf g_i(F_{-i})\neq \emptyset$. Further, for any $F_i\in \mathsf g_i(F_{-i})$, the matrices $A_i+ B_iF_i $ and $A_i+B_iF_i+G_iP_i$ are stable. Moreover, for any feasible $P_i\in \mathsf {S}_i$ and $M_i\in \mathsf {M}_i (P_i)$, an  $F_i\in \mathsf g_i(F_{-i})$ is obtained by solving the following LMIs
		\begin{align} 
			&\begin{bmatrix} (A_i+B_iF_i)^\prime P_i+P_i(A_i+B_iF_i)  + Q_i +P_i G_i P_i   &  F_i^\prime \\
				F_i & - R_i^{-1} \end{bmatrix} \prec 0,
			\label{eq:FGC2}\\
			&(A_i+B_iF_i+G_iP_i)^\prime M_i +M_i(A_i+B_iF_i+G_iP_i)\prec 0.
			\label{eq:stability2}
		\end{align} 
		\label{eq:OFstab22}
	\end{subequations} 
	The worst-case disturbance signal from Player $i$'s perspective  is given by $\bar{d}_i(t)=D_i^{-1}E^\prime P_i e^{(A_i+B_iF_i+G_iP_i)t}x_0$ for $t\in[0,\infty)$. 
	\label{thm:GCoi2}
	\label{thm:synthSFB_disturbance}
\end{corollary}
\begin{proof}
The proof is similar to that of Theorem \ref{thm:GCoi}, and we briefly mention the key steps.  Pre and post multiplying $P_i^{-1} \succ 0$ in the BMI \eqref{eq:SCSGCE1}, we get $
	P_i^{-1}A_i^\prime + A_i P_i^{-1} + P_i^{-1}F_i^\prime B_i^\prime + B_i F_i P_i^{-1} + P_i^{-1} Q_i P_i^{-1} + P_i^{-1} F_i^\prime R_i F_i P_i^{-1} + G_i \prec 0$. 
Using the Schur complement and notation provided in the theorem statement, the previous inequality is equivalently written as
$	\Phi_i^{1} (P_i^{-1}) +  (\bar{Y}_i)^\prime  F_i^\prime B_i^{3} + (B_i^{3})^\prime  F_i \bar{Y}_i \succ 0$
where $\bar{Y}_i = \begin{bmatrix} P_i^{-1} &0_{n\times n} &0_{n\times m_i}\end{bmatrix} $. From Lemma \ref{lem:Finsler} the previous BMI (in $(F_i,P_i^{-1})$)  is feasible if and only if the following LMI (in $P_i^{-1}$) 
\begin{align}
	\mathrm N_{B_i^3}^\prime {\Phi}^1_i(P_i^{-1}) \mathrm N_{B_i^3}   \prec 0 
	\label{eq:projectionLMI_Pi_inv}
\end{align}
is feasible. Using Schur complement, the inequality \eqref{eq:SCSGCE2} is written as 
\begin{align}
	\begin{bmatrix}
		\delta_i  &x_0^\prime \\
		x_0     &P_i^{-1}
	\end{bmatrix}
	\succ 0.
\end{align}
If $\mathsf{P}_i \neq \emptyset $, then from \eqref{eq:setP2}, there is a $P_i \succ 0$ such that $P_i^{-1} \in \mathsf{Y}_i$. Using Lemma \ref{lem:Inq_upperbound}, we write  $(A_i+B_iF_i+G_iP_i)^\prime M_i +M_i(A_i+B_iF_i+G_iP_i) \preceq 
	(A_i+B_iF_i)^\prime M_i + M_i(A_i+B_iF_i) + \gamma_{i} P_i G_i P_i   + \frac{1}{\gamma_{i}} M_i G_i M_i$ 
for some $\gamma_i >0$. So, the satisfaction of the stricter inequality $(A_i+B_iF_i)^\prime M_i + M_i(A_i+B_iF_i) +
	\gamma_{i} P_i G_i P_i   + \frac{1}{\gamma_{i}} M_i  G_i M_i \prec 0$
implies satisfaction of  \eqref{eq:SCSGCE3}. Pre and post multiplying $M_i^{-1} \succ 0$ the previous inequality we get 
$	M_i^{-1} A_i^\prime + A_i M_i^{-1} + M_i^{-1} F_i^\prime B_i ^\prime + B_iF_i M_i^{-1} + \gamma_i M_i^{-1} P_i G_i P_i M_i^{-1} + \frac{1}{\gamma_{i}} G_i \prec 0$. If the set $\mathsf{P}_i \neq \emptyset $, then for a feasible $P_i \in \mathsf{P}_i$, the previous BMI in the variables $(M_i^{-1},F_i)$ is rewritten using Schur complement and the notation provided in the theorem statement as $\Phi_i^{2}(M_i^{-1};P_i) +  (\bar{V}_i)^\prime F_i^\prime B_i^{4} +( B_i^{4})^\prime F_i \bar{V}_i \prec 0$, where $\bar{V}_i = \begin{bmatrix} M_i^{-1} &0_{n\times q} \end{bmatrix}$. Next, using Lemma \ref{lem:Finsler}, the previous BMI is feasible (in $(M_i^{-1},F_i)$) if and only if the following LMI is feasible (in $M_i^{-1}$)
\begin{align}
	\mathrm N_{B_i^4}^\prime{\Phi}^2_i(M_i^{-1};P_i) \mathrm N_{B_i^4}\prec 0. 
	\label{eq:projectionLMI_Mi_inv}
\end{align}
If $\mathsf{S}_i \neq \emptyset$, then there is a $P_i \succ 0$ such that $P_i \in \mathsf{P}_i$ and $\mathsf{M}_i(P_i) \neq \emptyset$, which further implies that there is a $M_i \succ 0$ such $M_i \in \mathsf{M}_i(P_i)$. From \eqref{eq:setM2}, this means, $M_i \succ 0$ satisfies $M_i^{-1} \in V_i(P_i)$. Then following Lemma \ref{lem:Finsler}, there exists a $F_i \in \mathbb{R}^{m_i \times n}$ which is feasible for \eqref{eq:projectionLMI_Pi_inv}  and \eqref{eq:projectionLMI_Mi_inv} for any $P_i \in \mathsf{S}_i$ and $M_i \in \mathsf{M}_i(P_i)$. 
\end{proof}

\begin{remark}
 In the output feedback case, although the sets $\mathsf{X}_i$, $\mathsf{Y}_i$, $\mathsf{U}_i(Z)$, and $\mathsf{V}_i(Z)$ characterized by LMIs are convex (see \eqref{eq:XLMI2}, \eqref{eq:YLMI2}, \eqref{eq:ULMI2}, and \eqref{eq:VLMI2}), the sets $\mathsf{P}_i$, $\mathsf{M}_i(Z)$, and $\mathsf{S}_i$ (see \eqref{eq:setP2_OF}, \eqref{eq:setM2_OF}, and \eqref{eq:setPs2_OF}) are observed to be non-convex due to the coupling constraints $P_i \in \mathsf{X}_i$, $P_i^{-1} \in \mathsf{Y}_i$, $M_i(Z) \in \mathsf{U}_i(Z)$, and $M_i^{-1}(Z) \in \mathsf{V}_i(Z)$. However, in the state feedback case, we note that these sets (see \eqref{eq:setP2}, \eqref{eq:setM2}, and \eqref{eq:setPs2}) are convex due to the absence of coupling constraints. This observation implies that SCSGCE strategies can be synthesized using LMIs without the need for SDP relaxation schemes as used in Section~\ref{sec:algorithm}.
 
\end{remark}

Again, for computational purposes, we prefer closed sets. So, we consider the following closed ``$\epsilon$-approximation" of the convex, open sets $\mathsf{Y}_i$ and $\mathsf{V}_i(Z)$ respectively as 
$\mathsf Y^\epsilon_i :=\big\{Y \in \mathbb R^{n\times n}~\big|~Y\succ 0,~\left[\begin{smallmatrix}\delta_i &x_0^\prime\\x_0&Y\end{smallmatrix}\right] \succeq  \epsilon I_{n+1},  ~ \mathrm N_{B_i^3}^\prime {\Phi}_i^1(Y)\mathrm N_{B_3^1}\preceq -\epsilon I_{2n+m_{i}} \big\}$ and $\mathsf V^\epsilon_i(Z):=\big\{V\in \mathbb R^{n\times n} ~\big|~V\succ 0,~\mathrm N_{B_i^4}^\prime \Phi_i^2(V;Z) \mathrm N_{B_i^4}\preceq -\epsilon I_{n+q} \big\}$ for some $\epsilon >0$.
Then, using this, the sets defined by \eqref{eq:setP2}-\eqref{eq:setPs2} are approximated as $\mathsf{P}^\epsilon_i:= \{P_i\in \mathbb{R}^{n\times n} ~ |~P_i\succ 0, P_i^{-1}\in \mathsf{Y}^\epsilon_i \}$, $\mathsf{M}^\epsilon_i(Z):= \{M_i \in \mathbb{R}^{n\times n} ~|~ M_i\succ 0,M_i^{-1}\in \mathsf{V}^\epsilon_i(Z) \}$, and $\mathsf{S}^\epsilon_i:= \{P_i\in \mathsf{P}^\epsilon_i~|~\mathsf{M}^\epsilon_i(P_i) \neq \emptyset \}$. To assess the feasibility of $\mathsf{g}_i(F_{-i}) \neq \emptyset$, for a given $(\delta_i,F_{-i})$, we introduce Algorithm \ref{alg:3}. To begin, the initial stage involves verifying if $\mathsf{P}^\epsilon_i\neq \emptyset$, which is an LMI feasibility problem. If this yields a solution, denoted as $P_i\in \mathsf{P}^\epsilon_i\neq \emptyset$, the next step involves verifying if $\mathsf{M}_i^\epsilon (P_i)\neq \emptyset$, which again is an LMI feasibility problem. If this yields a solution, then we conclude $\mathsf{S}_i^\epsilon \neq \emptyset$, and Player $i$'s feedback strategy, ensuring the worst-case cost $J_i^{\mathrm{sc}}(F_i,F_{-i})< \delta_i$, is synthesized from \eqref{eq:FGC2}-\eqref{eq:stability2}.
Next, to find an SCSGCE which satisfies the required sufficient conditions, we use sequential guaranteed cost response approach as outlined in Algorithm \ref{alg:2}.  {Specifically, in Step 2 of Algorithm \ref{alg:2}, we use Algorithm \ref{alg:3}. In Step 3, we obtain $F_i$ using \eqref{eq:FGC2}-\eqref{eq:stability2} for a feasible $P_i \in \mathsf{S}^\epsilon_i \neq \emptyset$.}

\begin{algorithm}[H]
	\captionsetup{labelfont={sc,bf}, labelsep=newline}
	\setstretch{1}
	\caption{SCSGCE: Non-emptiness of the guaranteed cost response} 	\label{alg:3}
	Determine $P_i \in \mathsf{P}_i^\epsilon$ 
	\\
	\eIf{ $\mathsf{P}_i^\epsilon \neq \emptyset $}
	{Determine $M_i \in \mathsf{M}_i^\epsilon(P_i)$ \\
		\eIf{ Found $(P_i,M_i)$ s.t. $P_i \in \mathsf{P}_i^\epsilon$ and $M_i \in \mathsf{M}^\epsilon_i(P_i)$, this implies $\mathsf{S}^\epsilon_i \neq \emptyset $}
		{$\mathsf{g}_i(F_{-i}) \neq \emptyset$}
		{Stop. $\mathsf{g}_i(F_{-i}) = \emptyset$}
	}
	{Stop. $\mathsf{g}_i(F_i) = \emptyset$}
\end{algorithm}
\section{Numerical illustrations}
\label{sec:numerical} 
In this section, we illustrate the performance of guaranteed cost equilibrium-based controls for the control of networked multi-agent systems.
\begin{example} \label{ex:exp2}
	\begin{figure}[h] \centering 
		\subfloat[SCOGCE] {\adjustbox{raise=-0.3pc}{\includegraphics[scale=.6,valign=t]{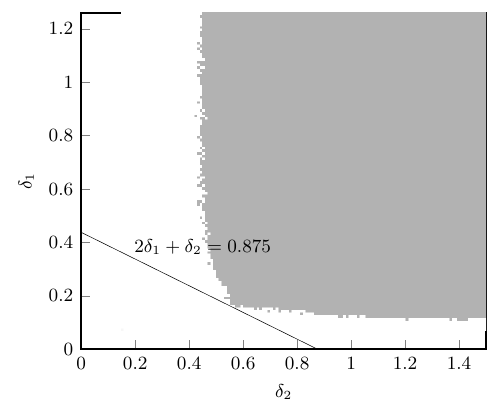}  \label{fig:fig2}}} \quad 
		\subfloat[PoS vs $\frac{2\delta_1+\delta_2}{J^\mathrm{co}}$]{\includegraphics[scale=.6,valign=t]{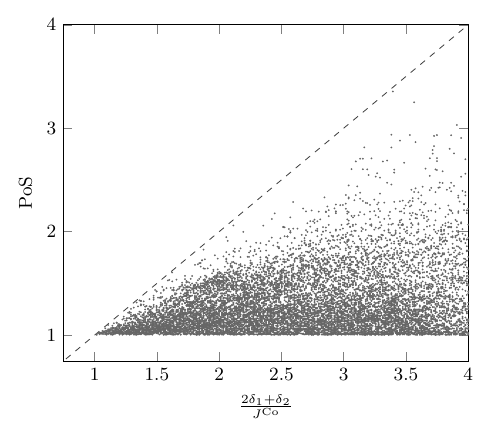}
			\label{fig:fig31}
		} \quad   
		\subfloat[PoS vs $\frac{2\delta_1+\delta_2}{J^\mathrm{co}}$]{\includegraphics[scale=.6,valign=t]{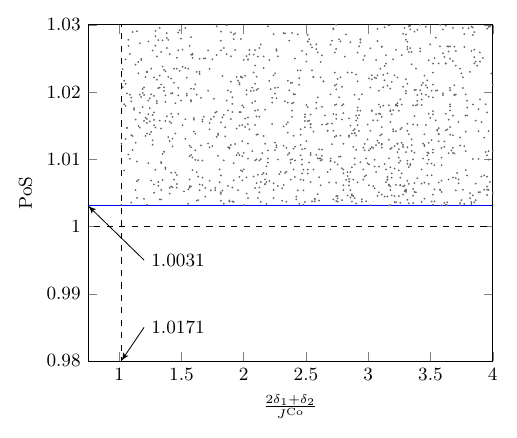}
			\label{fig:fig32}}     
		\caption{In Panel (a), the shaded region in the $\delta_1-\delta_2$ plane corresponds to the parameter region where an SCOGCE exists. Panel(c) is a magnified version of panel (b) close to PoS=1.} 
		\label{fig:fig3}
	\end{figure}

We consider the differential game presented in Example \ref{ex:exp1}. Due to symmetry, we assume $\delta_1 = \delta_3$. For numerical simulations,  {we consider the region $(\delta_1, \delta_2) \in [0.07, 1.25] \times [0.15, 1.45]$} and uniformly sample this region with a resolution of $0.01$. For each sample point in this grid, using Algorithms \ref{alg:1} and \ref{alg:2}, we check for the existence of SCOGCE. The shaded portion in Figure \ref{fig:fig2} corresponds to the region in the $\delta_1 - \delta_2$ plane where an SCOGCE exists. At an SCOGCE $(F_i^\circ, F_{-i}^\circ)$, the individual cost of players satisfies $J^{\mathrm{sc}}_i(F_i^\circ, F_{-i}^\circ) < \delta_i$ for $i = 1, 2, 3$. Then, from \eqref{eq:PoSdef}, this implies $\mathrm{PoS} \leq \frac{2\delta_1 + \delta_2}{J^\text{Co}}$. Using \eqref{eq:tgcostcoop}, the total game cost in cooperation is obtained as $J^\text{Co} = 0.875$. We note that the shaded region in Figure \ref{fig:fig2} lies above the line $2\delta_1 + \delta_2 = 0.875$, indicating that the upper bound on the Price of Stability (PoS) associated with SCOGCE satisfies $\text{PoS} \leq \frac{2\delta_1 + \delta_2}{0.875} > 1$. Figure \ref{fig:fig31} illustrates the relationship between PoS and the upper bound $\frac{2\delta_1 + \delta_2}{0.875}$ for SCOGCE, verifying the previous relation. Figure \ref{fig:fig32} illustrates that there do not exist SCOGCEs for the parameter values  {in the region $ (\delta_1, \delta_2) \in \{[0.07, 1.25] \times [0.15, 1.45]~|~\frac{2\delta_1 + \delta_2}{J^\text{Co}} < 1.0171 \}$.} Furthermore, Figure \ref{fig:fig32} also illustrates that the lowest PoS achievable by SCOGCE is $1.0031$. This implies that there exist many SCOGCEs that achieve welfare levels close to the one obtained in the cooperative outcome; see Remark \ref{rem:PoS}.
\end{example}  
\begin{example}
	\label{ex:exp4_FS_distb}
	\begin{figure}[H] 
		\centering 
		\begin{tikzpicture}[scale=.25,>=latex', inner sep=1mm, font=\small]
			\tikzstyle{solid node}=[circle,auto=center,draw,minimum size=6pt,inner sep=2,fill=black!6]
			\tikzstyle{dedge} = [draw, blue!50, -,  line width=.5mm]
			\node (n1)[solid node] at (-4,0) {\scriptsize{$1$}};
			\node (n2)[solid node] at (4,0) {\scriptsize{$2$}};
			\node (n3)[solid node] at (0,-4) {\scriptsize{$3$}}; 
			\path[dedge] (n1) --  (n2);
			\path[dedge] (n2) --  (n3);	 
			\path[dedge] (n1) --  (n3);
		\end{tikzpicture}
		\caption{3-agent networked multi-agent system with full state observations}
		\label{fig:3agentnetwork_FS_distb}
	\end{figure}
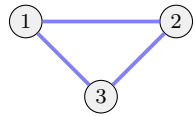
	
	We consider a 3-agent networked multi-agent system, as shown in Figure \ref{fig:3agentnetwork_FS_distb}, with complete state observations.  {The state dynamics associated with agent $i \in \mathsf{N} = \{1,2,3\}$ are given by 
		$\dot{x}_i(t)= a_{ii} x_i(t)+b_{ii} u_i(t)+ \sum_{j \in \mathsf{N}_i} b_{ij} u_j(t) + e_i d(t)$, 
		$a_{11} = -2$, $a_{22} = 2$, $a_{33} = 1$, $b_{11} = 1$, $b_{12} = 0.6$, $b_{13} = 0.3$, $b_{21} = 0.2$, $b_{22} = 0.8$, $b_{23} = 0.4$, $b_{31} = 0.5$, $b_{32}= 0.2$, $b_{33} = 0.6$, $e_1 = 0.4$, $e_{2} = 0.5$ and $e_{3} = 0.3$. Here, $x_i(t) \in \mathbb{R}$ and $u_i(t) \in \mathbb{R}$ denote the state and control signals of agent $i$. }The global state vector is given by $x(t) := \col{x_i(t)}, i \in \mathsf{N}$, and the output vector of agent $i \in \mathsf{N}$ is given by $y_i(t) := x(t)$. The objectives of the agents are given by $J_1 = \int_0^\infty ((x_1(t) - x_2(t))^2 + (x_2(t) - x_3(t))^2 + u_1^2(t) - 4d^2(t))~dt$, $J_2 = \int_0^\infty ((x_2(t) - x_1(t))^2 + (x_2(t) - x_3(t))^2 + u_2^2(t) - 4d^2(t))~dt$, and $J_3 = \int_0^\infty ((x_3(t) - x_2(t))^2 + (x_3(t) - x_1(t))^2 + u_3^2(t) - 4d^2(t))~ dt$. Upon writing the objectives and dynamics in standard form \eqref{eq:system}, we get $Q_{1}=\left[\begin{smallmatrix} 2 &-1 &-1 \\ -1 &1 &0 \\ -1 &0 &1\end{smallmatrix}\right]$ , $Q_{2} = \left[\begin{smallmatrix} 1 &-1 &0 \\ -1 &2 &-1 \\ 0 &-1 &1 \end{smallmatrix}\right]$, and
	$Q_{3}=\left[\begin{smallmatrix} 1 &0 &-1 \\ 0 &1  &-1 \\ -1 &-1 &2 \end{smallmatrix}\right]$.  For simulation purposes, we choose $x_0 = [0.2~~ -0.3~~ 0.4]^\prime$, and $R_i = 1$, $D_i = 4$ for $i \in \mathsf{N}$. Solutions of CARE involve solving $18$ multivariate polynomials of degree $2$ in $18$ variables. We obtain one stabilizing solution $(P_1, P_2, P_3)$ of CARE, that is, the solution of \eqref{eq:care1} such that the matrices $A-\sum_{j=1}^{3} B_jR_j^{-1} B_j^\prime P_j$ and $A- \sum_{j=1}^{3} B_j R_{j}^{-1} B_j^\prime P_j + ED_i^{-1}E^\prime P_i$ for all $i \in \mathsf{N}$ are stable. The soft-constrained state feedback Nash equilibrium (SCFNE) costs of the agents are obtained as $(J_1^{\star}, J_2^{\star}, J_3^{\star}) = (0.0999, 2.3385, 1.9535)$, and the total game cost at SCFNE is obtained as $J^{*} = \sum_{i=1}^{3} J_i^{\star} = 4.3918$. The total game cost in cooperation is obtained as $J^{\text{Co}} = 3.2801$ by solving the problem \eqref{eq:tgcostcoop}. So, the $\mathrm{PoS}$ associated with the SCFNE is computed as $\mathrm{PoS} = \tfrac{J^{*}}{J^{\text{Co}}} = 1.3390$. We assume $\delta_i = \delta$ for $i \in \mathsf{N}$ and sample the region $\delta \in [1.1,~6]$ with a resolution of $0.01$. For each sample point in this grid, using Algorithms \ref{alg:2} and \ref{alg:3}, we check for the existence of SCSGCE. Figure \ref{fig:posplot_3player_SF1_distb} illustrates the plot of $\mathrm{PoS}$ at SCSGE versus $\delta$, and this indicates that there does not exist an SCSGCE for $\delta < 1.85$. Figure \ref{fig:posplot_3player_SF2_distb} illustrates that there exist many SCSGCEs at which $\mathrm{PoS}$ is lower than in comparison with the $\mathrm{PoS}$ at SCFNE (denoted by the red line). The minimum achievable $\mathrm{PoS}$ at an SCSFGCE is computed as $1.0208$ (denoted by the blue line in Figure  \ref{fig:posplot_3player_SF2_distb}). This observation suggests that there are several SCSGCEs that attain welfare levels close to the one obtained with the cooperative outcome; see Remark \ref{rem:PoS}.
	\begin{figure}[h] \centering 
		\subfloat[$\mathrm{PoS}$ vs $\delta$]
		{\input{posplot_3player_SF1_distb} \label{fig:posplot_3player_SF1_distb}} \quad \quad  
		\subfloat[$\mathrm{PoS}$ vs $\delta$]{\input{posplot_3player_SF2_distb}
			\label{fig:posplot_3player_SF2_distb}
		} 
		\caption{Panel (a) illustrates $\mathrm{PoS}$ vs $\delta$ plot, and Panel (b) is a magnified version of Panel (a) close to $\mathrm{PoS} = 1$} 
		\label{fig:pos_ex4}
	\end{figure}
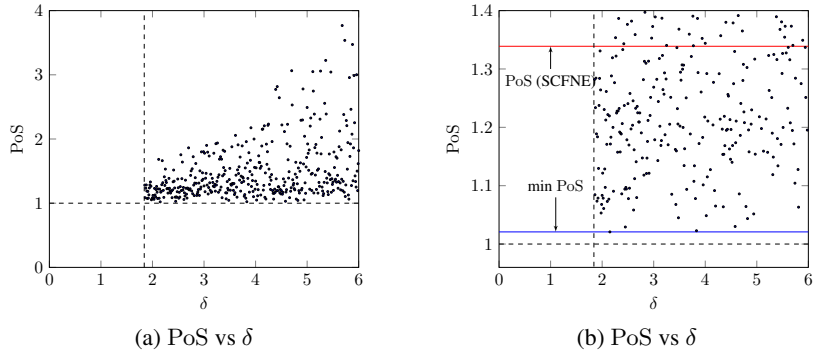
\end{example}
\begin{example}
		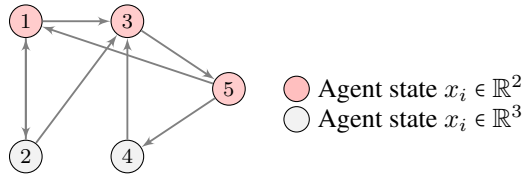
\begin{figure}[h] 
		\centering 
		\begin{tikzpicture}[scale=.225,>=latex', inner sep=1mm, font=\small]
			\tikzstyle{solid node1}=[circle,auto=center,draw,minimum size=6pt,inner sep=2,fill=red!20]
			\tikzstyle{solid node2}=[circle,auto=center,draw,minimum size=6pt,inner sep=2,fill=black!5]
			\tikzstyle{dedge} = [draw, ->, black!50,  line width=.25mm]
			\tikzstyle{leader node}=[circle,auto=center,draw,minimum size=6pt,inner sep=2,fill=blue!25]
			\def \x{2.5};
			\def \y{4};
			
			\node (n1)[solid node1] at (6,\y) {\scriptsize{$1$}};
			\node (n2)[solid node2] at (6,-\y) {\scriptsize{$2$}}; 
			\node (n3)[solid node1] at (12,\y) {\scriptsize{$3$}};
			\node (n4)[solid node2] at (12,-\y) {\scriptsize{$4$}};
			\node (n5)[solid node1] at (18,0) {\scriptsize{$5$}};  
			\path[dedge] (n1)  --  (n2);
			\path[dedge] (n2) --  (n1);	 
			\path[dedge] (n1) --  (n3);	 
			\path[dedge] (n5) --  (n4);	 
			\path[dedge] (n3) --  (n5);	 
			\path[dedge] (n4) --  (n3);	
			\path[dedge] (n2) --  (n3);	
			\path[dedge] (n5) --  (n1);	
			\draw[fill=red!25](22,0) circle[radius=0.75];
			\draw[fill=black!5](22,-1.8) circle[radius=0.75];
			\draw(22,0)node[right,xshift=5]{Agent state $x_i \in \mathbb{R}^2$};
			\draw(22,-1.8)node[right,xshift=5]{Agent state $x_i \in \mathbb{R}^{3}$};
		\end{tikzpicture}
		\caption{5-agent networked multi-agent system.}
		\label{fig:5agentnetwork}
	\end{figure}
	In this example, we demonstrate the performance of SCOGCE controllers for output consensus of heterogeneous networked multi-agent systems. We consider a 5-agent networked multi-agent system where the agents are connected according to the directed graph as illustrated in Figure \ref{fig:5agentnetwork}. The in-neighbor sets of the agents are  $\mathsf{N}_1 = \{2,5\}, \mathsf{N}_2 = \{1\}, \mathsf{N}_3 = \{1,2,4\}, \mathsf{N}_4 = \{5\}$, and $\mathsf{N}_5 = \{3\}$. So, the in-degree of each player is $d_{1} = 2$, $d_{2} = 1$, $d_{3} = 3$, $d_{4} = 1$ and $d_{5} = 1$. The  dynamics of agent $i$ is affected by an external deterministic disturbance signal, and is given by $\dot{x}_i(t)=\bar{A}_ix_i(t)+\bar{B}_iu_i(t)+\bar{E}_id(t)$, $\bar{y}_i(t)=\bar{C}_ix_i(t)$, 
	\begin{align*}
		&\bar{A}_{i} = \begin{bmatrix}	0.3 &-2 \\ 0.1 &-0.2	 \end{bmatrix},~
		\bar{B}_{i} = \begin{bmatrix} 	1.8 &-0.8 \\ 0.9 &1.6 \end{bmatrix},
		\bar{C}_{i} = \begin{bmatrix} 	-0.1 &1.2 \\ 0.4 &1.4 \end{bmatrix}, \bar{E}_{i} = \begin{bmatrix} 1 \\ 0\end{bmatrix}~i=1,3,5,\\
		&\bar{A}_{i} =\begin{bmatrix}  	0 &1 &0 \\ 0 &0 &1 \\ 0 &0 &-2\end{bmatrix},~
		\bar{B}_i=\begin{bmatrix}     6 &0 \\ 0 &1 \\ 1 &0 \end{bmatrix},~
		\bar{C}_i=\begin{bmatrix}      1 &0 &0 \\ 0 &1 &0 \end{bmatrix},\bar{E}_{i} =  \begin{bmatrix} 0.5 \\-1 \\1 \end{bmatrix}~i=2,4,
	\end{align*}
	where $x_{i}(t)$,  $u_i(t)$ and $\bar{y}_i(t) $ denote respectively the state, control and output vectors of agent $i$. 
	The global state vector is given by $x(t):=\col{x_i(t)}_{i\in \mathsf N}$, and the local output information of agent $i\in \mathsf N$ is $y_i(t) := \col{\bar{y}_k(t)}_{k \in \mathsf{N}_i \cup i} = C_i x(t)$, where $C_i \in \mathbb{R}^{ s_{i} \times n}$. The local reference is given by $r_i(t)=\frac{1}{ d_i} \sum_{j\in \mathsf N_i} \bar{y}_j(t) $. The difference between the output and the local reference signal is given by $e_i(t)=\bar{y}_i(t)-r_i(t)$. Now, we write the error signal $e_i(t)$ for each agent as $e_1(t) =\bar{y}_1(t)-\tfrac{1}{2}\bar{y}_2(t)-\tfrac{1}{2}\bar{y}_5(t)=\left(\begin{bmatrix}1 &-\tfrac{1}{2} & -\tfrac{1}{2} \end{bmatrix} \otimes I_2 \right)y_1(t)$, $e_2(t) = \bar{y}_2- \bar{y}_1= \left(\begin{bmatrix} -1 &1\end{bmatrix} \otimes I_2\right)y_2(t)$, $e_3(t) = \bar{y}_3-\tfrac{1}{3}\bar{y}_1-\tfrac{1}{3}\bar{y}_2-\tfrac{1}{3}\bar{y}_4 = \left(\begin{bmatrix} -\tfrac{1}{3} &-\tfrac{1}{3} &1 &-\tfrac{1}{3}\end{bmatrix} \otimes I_2 \right)y_3(t) $, $e_4(t) = \bar{y}_4(t) - \bar{y}_5(t) = \left(\begin{bmatrix}
		1 &-1 \end{bmatrix} \otimes I_2\right)y_4(t)$, and $e_5(t) = \bar{y}_5(t)-\bar{y}_3(t) = \left(\begin{bmatrix} -1 &1 \end{bmatrix} \otimes I_2\right)y_5(t)$. The objective of  each agent $i \in \mathsf{N}$ is given as follows 
	\begin{align*}
		J_i  := \int_{0}^{\infty} \left(e_i^\prime(t) e_i(t) + u_i^\prime(t) R_i u_i(t) - d^\prime(t) D_i d(t) \right)~dt. \label{eq:objlocaloutputEx4}
	\end{align*}
	The control strategy of each agent $i \in \mathsf{N}$ in the above objective is to drive the error to zero, that is, $ e_i(t) = 0$ as $t \rightarrow \infty, \forall i \in \mathsf{N}$. The objectives are rewritten in the standard form \eqref{eq:system} as follows 
	\begin{align*}
		J_i(u_i,u_{-i},d) := \int_{0}^{\infty} \left(y_i^\prime(t) Q_i y_i(t) + u_i ^\prime(t) R_i u_i(t) - d^\prime(t) D_i d(t) \right)~dt,
	\end{align*}
	where $Q_{1} = \left[\begin{smallmatrix} 	1& -\frac{1}{2} & -\frac{1}{2}\\
		-\frac{1}{2} & \frac{1}{4} & \frac{1}{4} \\	-\frac{1}{2} &  \frac{1}{4} &  \frac{1}{4} 	\end{smallmatrix}\right]\otimes I_2$, 
	$Q_2 =  \left[\begin{smallmatrix}1&-1\\-1&1  \end{smallmatrix} \right]  \otimes I_2$, 
	$Q_3 =   \left[\begin{smallmatrix}\frac{1}{9} &\frac{1}{9} &-\frac{1}{3} &\tfrac{1}{9}\\
		\frac{1}{9} &\frac{1}{9} &-\frac{1}{3} &\frac{1}{9}\\
		-\frac{1}{3} &-\frac{1}{3} &1 &-\frac{1}{3}\\
		\frac{1}{9} &\frac{1}{9} &-\frac{1}{3} &\frac{1}{9} \end{smallmatrix}\right]  \otimes I_2$, 
	$Q_4 =  \left[ \begin{smallmatrix}1&-1\\-1&1  \end{smallmatrix}\right]  \otimes I_2$, $Q_5 = \left[\begin{smallmatrix}1&-1\\-1&1  \end{smallmatrix} \right] \otimes I_2$.  For simulation purpose, we choose $(\delta_1,\delta_2,\delta_3,\delta_4,\delta_5)=(1.30, 0.7, 1.30, 0.7, 1.30)$, 
	 $x_0=[-0.3~~0.5~~0.4~~0.2~~0.6~~-0.3~~0.2~~-0.1~~ 0.5~~ 0.7~ ~0.2~~ -0.4]^\prime$, $R_i=I_2$, $i\in \mathsf N$, and  $D_1=D_3=D_5=12.5$ and $D_2=D_4=10$. We assume disturbance signal is $d(t) = 10\sin{t}e^{-t} \in L_2([0,\infty))$. Using the  Algorithms \ref{alg:1} and \ref{alg:2}, we found the following set of SCOGCE strategies 
	\begin{align*}
		&F_{1}^{\circ} = 
		\begin{bmatrix}
			7.4397  &-5.4180   &-1.8863   &-0.2325    &5.6512   &-6.1497 \\
			-2.7672   &-1.3434   &0.6394    &1.1094   &-0.8136   &-0.9275 
		\end{bmatrix}, \\
		&F_2^{\circ} = 
		\begin{bmatrix}
			0.7508    &0.9721   &-5.0861   &-0.0788 \\
			2.2360   &-0.2057  &-2.7248   &-3.6266
		\end{bmatrix}, \\
		&F_3^{\circ} = 
		\begin{bmatrix}
			-0.7364  & -0.0965   & 0.1044   &-0.1584    &7.3679   &-5.6129   & 1.3109   &-0.3174 \\
			1.3736   &-0.6591    &0.3946    &1.2000   &-3.3657   &-7.1549   &-1.1549    &1.6549 
		\end{bmatrix}, \\
		&F_4^{\circ} = 
		\begin{bmatrix}
			-3.6049  & -1.5738  & -0.9251    &0.6961 \\
			3.0475   &-1.7145   &-8.8981    &9.2294 
		\end{bmatrix}, \\
		&F_5^{\circ} = 
		\begin{bmatrix}
			5.6955   &-4.0258   &17.6970   &-17.2408 \\
			0.2054   &-0.5258   &-0.4322   &-3.9581
		\end{bmatrix}.
	\end{align*}
	Agents' worst-case cost at SCOGCE $(J^{\mathrm{sc}\circ}_1,J^{\mathrm{sc}\circ}_2,J^{\mathrm{sc}\circ}_3,J^{\mathrm{sc}\circ}_4,J^{\mathrm{sc}\circ}_5) = (0.9848, 0.0547, 0.1113, 0.6687, 0.6895)$, and the total game cost at SCOGCE is computed as $ 2.5090$. Next, we compute the total cost in cooperation which is $J^{\text{Co}} = 1.9993$ (see \eqref{eq:tgcostcoop}), which implies the  $\mathrm{PoS} $ at SCOGCE is  $1.2549$.  We observe that in Figure \ref{fig:ei_1_plot_case1magn_ex3}, and Figure \ref{fig:ei_2_plot_case1magn_ex3} error trajectories of the agents converge to zero asymptotically. Figure \ref{fig:disturbance_case1_ex3} illustrates the worst-case disturbance trajectories from the perspectives of players $i$. 
	
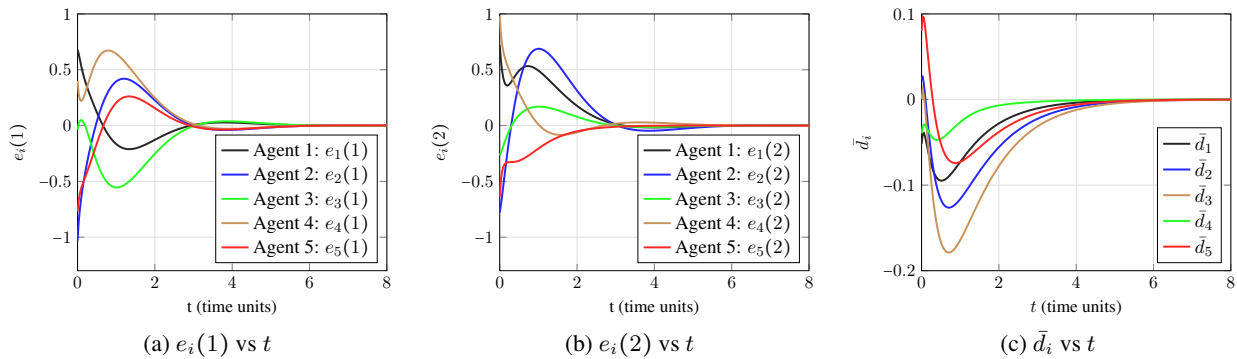
\begin{figure}[h] \centering 
	\subfloat[$e_i(1)$ vs $t$]{\input{ei_1_plot_case1magn_ex3} \label{fig:ei_1_plot_case1magn_ex3}}~  
	\subfloat[$e_i(2)$ vs $t$]{\input{ei_2_plot_case1magn_ex3} \label{fig:ei_2_plot_case1magn_ex3}}~   
	\subfloat[$\bar{d}_i$ vs $t$]{\input{disturbance_case1_ex3}\label{fig:disturbance_case1_ex3}}     
	\caption{In Panel (a) and Panel (b) error trajectories of the agents using SCOGCE strategies. Panel (c) worst case disturbance trajectories from each player perspectives} 
	\label{fig:trajectories}
\end{figure}
\end{example}

\section{Conclusions}
\label{sec:conclusions}
In this paper, we studied a solution concept referred to as the soft-constrained output feedback guaranteed cost equilibrium for infinite-horizon linear-quadratic uncertain differential games with an output feedback information structure. These equilibrium strategies ensure that players' worst-case costs are upper-bounded by a given threshold. We provide sufficient conditions for the existence of these equilibria and also provide linear matrix inequality-based iterative schemes for synthesizing them.
For future work, it would be interesting to investigate the existence of guaranteed cost equilibria when the objectives of the players are hard-bound constrained, that is, when the disturbance signal has a finite norm.
\bibliographystyle{elsarticle-num} 
\bibliography{SCNERef2}
\appendix
\section{Appendix} 
In this appendix, we recall few results from linear systems theory and linear matrix inequalities.
\begin{lemma}[{\cite[Theorem 6.4.1]{Mareels:96}}]
	Consider the system
	\begin{align}
		\dot{x}(t)=Ax(t)+Bu(t),~x(0)=x_0,
	\end{align}
   where $x(t)\in \mathbb R^n$, $u(t)\in \mathbb R^m$, $A\in \mathbb R^{n\times n}$, $B\in \mathbb R^{n\times m}$, and $x_0\in \mathbb R^n$. 
	Assume that $A$ is stable. 
	\begin{enumerate}
		\item There exists constants $C_{0}$, $C_{1} \geq 0$ such that for all input functions $u$, all initial states $x_{0}$, and all $t_{0} \leq t$ the following holds
		\begin{align}
			\int_{t_{0}}^{t} x^\prime (\tau) x(\tau) d\tau \leq C_{0} x_{0}^\prime x_{0} + C_1 \int_{t_{0}}^{t} u^\prime (\tau) u(\tau) d\tau.
		\end{align}
		\item If $u (.)\in L^m_{2}([0,\infty))$, then $\lim_{t \rightarrow \infty} x(t) = 0 $. 
	\end{enumerate}
	\label{lem:L2stability}
\end{lemma}
\begin{lemma}[{\cite[Lemma 2.1]{Duan:13}}]
	Let $ {X}, {Y} \in \mathbb{R}^{m \times n}$, $L\in \mathbb R^{m\times m}$, ${L}\succ 0$, and $\gamma > 0$ be a scalar, then
	\begin{align}
		{X}^\prime  {L}  {Y} +  {Y}^\prime  {L} {X} \preceq \gamma  {X}^\prime  {L}  {X} + \frac{1}{\gamma}  {Y}^\prime  {L}  {Y}. 
	\end{align}
	\label{lem:Inq_upperbound}
\end{lemma}
\vskip25ex 
\end{document}

%% file: posplot_3player_SF1_distb.tex
\begin{tikzpicture}[scale=.6,>=latex']

\begin{axis}[
    xmin=0, xmax=6,
    ymin=0, ymax=4,
    xlabel={$\delta$},
    ylabel={$\mathrm{PoS}$},
    only marks,
    mark=*,
    mark size=0.6pt, 
    mark options={solid, black}
  ]
  \addplot table {dataEx4/pos_3player_SF_distb.dat};
  \draw[dashed,black, line width=0.2pt] (axis cs:0,1) -- (axis cs:6,1);
  \draw[dashed,black, line width=0.4pt] (axis cs:1.84,0) -- (axis cs:1.84, 4);
  \end{axis}
\end{tikzpicture}

%% file: posplot_3player_SF2_distb.tex
\begin{tikzpicture}[scale=.6,>=latex']

\begin{axis}[
    xmin=0, xmax=6,
    ymin=0.96, ymax= 1.4,
    xlabel={$\delta$},
    ylabel={$\mathrm{PoS}$},
    only marks,
    mark=*,
    mark size=0.6pt, 
    mark options={solid, black}
  ]
 \addplot table {dataEx4/pos_3player_SF_distb.dat};
  
  \draw[dashed,black, line width=0.4pt] (axis cs:0,1) -- (axis cs:6,1);
   \draw[solid,red, line width=0.4pt] (axis cs:0,1.3389) -- (axis cs:6,1.3389);
   \draw[solid,blue, line width=0.4pt] (axis cs:0,1.0208) -- (axis cs:6,1.0208);
   \draw[dashed,black, line width=0.4pt] (axis cs:1.84,0) -- (axis cs:1.84, 1.4);
   \draw[<-,black] (axis cs:1.1,1.0208) to (axis cs:1.1, 1.08) node[above]{{\small min $\mathrm{PoS}$}};
   \draw[<-,black] (axis cs:1,1.339) to (axis cs:1, 1.3) node[below]{{\small $\mathrm{PoS}$ (SCFNE)}};
  \end{axis}
\end{tikzpicture}

%% file: ei_1_plot_case1magn_ex3.tex
\begin{tikzpicture}[scale=.6,>=latex']
	\tikzset{every pin/.append style={font=\large}}
	\begin{axis}[xmin=0,xmax=8,ymin=-1.3, ymax=1, xlabel = {t (time units)},ylabel={$e_{i}(1)$},legend pos = south east,
		grid=both, legend style={nodes={scale=1.1 }}, 
		grid style={line width=.1pt, draw=gray!25},] 
		\addplot[color=black!85, line width=1.25pt] table{dataEx3/e1_1_ex3_SCOGCE.dat}; \addlegendentry{Agent 1: $e_{1}(1)$};
		\addplot[color=blue!85, line width=1.25pt] table{dataEx3/e2_1_ex3_SCOGCE.dat}; \addlegendentry{Agent 2: $e_{2}(1)$};
		\addplot[color=green!85, line width=1.25pt] table{dataEx3/e3_1_ex3_SCOGCE.dat}; \addlegendentry{Agent 3: $e_{3}(1)$};
		\addplot[color=brown!85, line width=1.25pt] table{dataEx3/e4_1_ex3_SCOGCE.dat}; \addlegendentry{Agent 4: $e_{4}(1)$}
		\addplot[color=red!85, line width=1.25pt] table{dataEx3/e5_1_ex3_SCOGCE.dat}; \addlegendentry{Agent 5: $e_{5}(1)$}
	\end{axis}  
\end{tikzpicture} 

%% file: ei_2_plot_case1magn_ex3.tex
\begin{tikzpicture}[scale=.6,>=latex']
	\tikzset{every pin/.append style={font=\large}}
	\begin{axis}[xmin=0,xmax=8,ymin=-1.3, ymax=1, xlabel = {t (time units)},ylabel={$e_{i}(2)$},legend pos = south east,
		grid=both, legend style={nodes={scale=1.1 }}, 
		grid style={line width=.1pt, draw=gray!25},] 
		\addplot[color=black!85, line width=1.25pt] table{dataEx3/e1_2_ex3_SCOGCE.dat}; \addlegendentry{Agent 1: $e_{1}(2)$};
		\addplot[color=blue!85, line width=1.25pt] table{dataEx3/e2_2_ex3_SCOGCE.dat}; \addlegendentry{Agent 2: $e_{2}(2)$};
		\addplot[color=green!85, line width=1.25pt] table{dataEx3/e3_2_ex3_SCOGCE.dat}; \addlegendentry{Agent 3: $e_{3}(2)$};
		\addplot[color=brown!85, line width=1.25pt] table{dataEx3/e4_2_ex3_SCOGCE.dat}; \addlegendentry{Agent 4: $e_{4}(2)$}
		\addplot[color=red!85, line width=1.25pt] table{dataEx3/e5_2_ex3_SCOGCE.dat}; \addlegendentry{Agent 5: $e_{5}(2)$}
	\end{axis}  
\end{tikzpicture} 

%% file: disturbance_case1_ex3.tex
\begin{tikzpicture}[scale=.6,>=latex']
	\tikzset{every pin/.append style={font=\large}}
	\begin{axis}[xmin=0,xmax=8,ymin=-0.2, ymax=0.1, xlabel = {$t$ (time units)},ylabel={$\bar{d}_{i}$},legend pos = south east,
		grid=both, legend style={nodes={scale=1.1 }}, 
		grid style={line width=.1pt, draw=gray!25},] 
		\addplot[color=black!85, line width=1.25pt] table{dataEx3/w1_case1_ex3.dat}; \addlegendentry{$\bar{d}_{1}$};
		\addplot[color=blue!85, line width=1.25pt] table{dataEx3/w2_case1_ex3.dat}; \addlegendentry{$\bar{d}_{2}$};
		\addplot[color=brown!85, line width=1.25pt] table{dataEx3/w3_case1_ex3.dat}; \addlegendentry{$\bar{d}_{3}$};
		\addplot[color=green!85, line width=1.25pt] table{dataEx3/w4_case1_ex3.dat}; \addlegendentry{$\bar{d}_{4}$}
		\addplot[color=red!85, line width=1.25pt] table{dataEx3/w5_case1_ex3.dat}; \addlegendentry{$\bar{d}_{5}$}
	\end{axis}  
\end{tikzpicture}